\documentclass[a4paper,reqno]{amsart} 
\usepackage{pdfpages}
\usepackage{amssymb}
\usepackage[mathcal]{euscript} 
\usepackage[cmtip,all]{xy}
\usepackage{color}
\usepackage{tikz-cd} 
\usepackage{ragged2e}

\newcommand{\bydef}{:=}

\newcommand{\tripleprod}{\{\cdot,\cdot,\cdot\}}

\newcommand{\id}{\mathrm{id}}
\newcommand{\ex}{\mathrm{ex}}

\newcommand{\sgn}{\mathrm{sgn}}

\newcommand{\cA}{\mathcal{A}}
\newcommand{\cB}{\mathcal{B}} 
\newcommand{\cC}{\mathcal{C}}
\newcommand{\cD}{\mathcal{D}}

\newcommand{\cM}{\mathcal{M}}

\newcommand{\cW}{\mathcal{W}}
\newcommand{\cZ}{\mathcal{Z}}

\newcommand{\NN}{\mathbb{N}} 
\newcommand{\ZZ}{\mathbb{Z}}

\newcommand{\FF}{\mathbb{F}} 

\DeclareMathOperator{\Hom}{\mathrm{Hom}}
\DeclareMathOperator{\End}{\mathrm{End}}

\DeclareMathOperator{\Aut}{\mathrm{Aut}}

\DeclareMathOperator{\AAut}{\mathbf{Aut}}

\DeclareMathOperator{\supp}{\mathrm{Supp}}

\DeclareMathOperator{\Ind}{\mathrm{Ind}}

\newcommand{\GLs}{\mathbf{GL}}

\newcommand{\subo}{\bar 0} 
\newcommand{\subuno}{\bar 1}

\newcommand{\invol}{\,-\,}

\newcommand{\Dinv}{D_{\mathrm{inv}}}

\newcommand{\Sgn}{\mathrm{Sgn}}

\newtheorem{theorem}{Theorem}[section]
\newtheorem{proposition}[theorem]{Proposition}
\newtheorem{lemma}[theorem]{Lemma}
\newtheorem{corollary}[theorem]{Corollary}

\theoremstyle{definition} 
\newtheorem{df}[theorem]{Definition}
\newtheorem{example}[theorem]{Example}
\newtheorem{examples}[theorem]{Examples}

\theoremstyle{remark} 
\newtheorem{remark}[theorem]{Remark}
\numberwithin{equation}{section}

\begin{document}
	
	\title{Gradings on associative triple systems of the second kind.}
	
	\author[A.~Daza-Garcia]{Alberto
		~Daza-Garcia} 
	\address{Departamento de
		Matem\'{a}ticas aplicadas (I), Universidad de Sevilla, 41012 Sevilla, Spain}
	\email{adaza1@us.es} 

	
	
	
	\begin{abstract} On this work we study associative triple systems of the second kind. We show that for simple triple systems the automorphism group scheme is isomorphic to the automorphism group scheme of the $3$-graded associative algebra with involution constructed by Loos. This result will allow us to prove our main result which is a complete classification up to isomorphism of the gradings of structurable algebras.
	\end{abstract}
	
	
	\maketitle

	\section{Introduction}\label{se:intro}

	Associative triple systems of the second kind (AT2), is a nonassociative structure which is reladed with other structures like  associative algebras with involution \cite{loos}, Lie and Jordan triple systems \cite{MeyLectures} or structurable algebras, concretely, with structurable algebras related to a hermitian form \cite{FauStr}. A complete description of AT2 was given by Otmar Loos \cite{loos} in 1971, who proved that all of them could be obtain as a subtriple system of the triple obtained in a concrete way from an associative algebra with involution $(A,\varphi)$ and a product given by $\{a,b,c\}=a\varphi(b)c$.
	
	The purpose of this work is to classify up to isomorphism the gradings by abelian groups on the simple finite-dimensional AT2. We will structure this paper into $6$ sections:
	
	In section \ref{sec:Prelim} we  generalize some of the tools given in \cite{EKmon} to work with gradings. Thus, we define $\FF$-linear $\Omega$-algebras as a generalization of linear nonassociative algebras which include for example triple systems or graded algebras with involution, and we give Proposition \ref{p:TransferThm}, which is a generalization of the classical transfer theorem for algebras.
	
	In section \ref{se:AT2} we define associative triple systems of the second kind, we recall the construction of O. Loos of an associative algebra with a $3$-grading and an involution and we give the construction of an AT2 from the associative algebras with Peirce decomposition. The main result of this section is the proof of existence of an isomorphism from the automorphism group scheme of a simple AT2 and the automorphism group scheme of the correspondent  associative algebra with a $3$-grading and an involution (Theorem \ref{t:AutSchemeIso}).
	
	In section \ref{sec:ArbitraryFields} we study the problem of classifying simple associative finite-dimensional triple systems over an arbitrary field of characteristic other than $2$.
	
	In section \ref{sec:AlgClosedField} we first classify graded division algebras with involution which are simple as algebras with involution over an algebraically closed field, (Theorems  \ref{t:DivAlgClass} and \ref{t:NonSimpleGradedDivAlg}) and finally we specialize the results of the previous section in Theorem \ref{t:ClassificationAssociativeCase}
	
	Finally, in section \ref{sec:GradingsAt2} we prove give the classification of simple AT2 in \ref{t:ClassificationTripleSystem}.
	
	On this paper, every vector space is assumed to be finite dimensional. We will assume $0\in \mathbb{N}$. The notation $G$ will stand for finitely generated abelian group and $e$ will be the neutral element.	$R$ will always denote an associative commutative unital ring, such that $\frac{1}{2}\in R$, $\FF$ will always be a field of characteristic different from $2$, and finally,  $\cM_n(R)$ will denote $n\times n$ matrices over $R$.

	\section{Preliminaries}\label{sec:Prelim}
	\subsection{$\Omega$-algebras}
	
	In this section $\cC$ will denote a monoidal category (see \cite{EGNOBook} for an introduction to monoidal categories). The reader who is not familiar with monoidal categories can assume $\cC$ to be the category of bimodules over an associative commutative unital ring $R$  with the normal tensor product and $\bold{1}=R$ (in which case $A^{\otimes 0}=R$).  In the case where $R=\FF$, we have the category of vector spaces.
	
	\begin{df}
		An \emph{operator domain} is a set $\Omega$ with a map $a\colon \Omega\to \NN$. The elements of $\Omega$ are called operators and for every $\omega\in \Omega$, $a(\omega)$ is called the \emph{arity} of $\Omega$ (or we say that $\omega$ is $a(\omega)$-\emph{ary}). We denote by $\Omega(n)$ the set of all the $n$-ary operators in $\Omega$.
	\end{df}
	
	\begin{df}
		Let $A$ be an object of $\cC$. An $\Omega$\emph{-algebra structure} on $A$ is a family of maps:
		
		\[\Omega(n)\to \Hom_{\cC}(A^{\otimes n},A).\] 
		We denote it by $(A,\Omega)$ or sometimes just $A$. If $\cC$ is the category of bimodules over an associative commutative ring $R$, we are going to call the $\Omega$-algebra structure $R$-\emph{linear} $\Omega$-\emph{algebra} (or just $R$-\emph{linear algebras} if $\Omega$ is known). Usually, for $\omega\in\Omega$ we are going to denote by $\omega$ its image under the previous map Given two $\Omega$-algebras $A$ and $B$ and $\omega\in \Omega$, we might denote by  $\omega^A$ and  $\omega^B$ the image of the operator for each respective algebras.
	\end{df}
	
	\begin{remark}If $\cC$ is the category of sets, the reader can learn about the theory of $\Omega$-algebras in \cite{CohUniversal}.
	\end{remark}
	
	\begin{example}\label{ex:ExamplesOmega} 
		\begin{itemize}
			\item[(1)] A nonassociative algebra is just an $\FF$-linear $\Omega$-algebra consisting on a vector space $A$ with a $2$-ary operator $\omega$. 
			\item[(2)] A ternary algebra $A$ is an $\FF$-linear $\Omega$-algebra consisting on a vector space $A$ with a $3$-ary operator.
			\item[(3)] An algebra with involution is an $\FF$-linear $\Omega$-algebra consisting on a vector space $A$ with a $2$-ary operator $\omega$ (the product) and a $1$-ary operator $\tau$, the involution, satisfying:
			\begin{itemize}
				\item[(i)] $\tau^2(x)=x$ for all $x\in A$
				\item[(ii)] $\tau(\omega(x,y))=\omega(\tau(y),\tau(x))$ for all $x,y\in A$.
			\end{itemize}   
			\item[(4)]
			Given an $\FF$-linear $\Omega$-algebra $A$, with a vector space decomposition  \[\Gamma\colon A=\bigoplus_{s\in S}A_s\]
			indexed by a set $S$, we can denote by $\Omega^{S}$, the operator domain in which \[\Omega^{S}(n)=\Omega(n)\] if $n\neq 1$ and \[\Omega^{S}(1)=\Omega(1)+\{\pi_s\mid s\in S\}.\]
			We can define an $\Omega^S$-algebra structure on $A$ by setting $\pi_s$ to be the projection on $A_s$ with respect to the vector space decomposition. 
			
		\end{itemize}
	\end{example}
	
	Given two $R$-linear $\Omega$-algebras $A$ and $B$, we say that $B$ is an $\Omega$-subalgebra of $A$ if $B\subseteq A$ and for any operator $\omega\in \Omega$, $\omega^A_{\mid B^n}=\omega^B$.
	
	Given two $R$-linear $\Omega$-algebras $A$ and $B$, an $R$-linear map $f\colon A\to B$ is a map of $\Omega$-algebras if for every $\omega \in \Omega(0)$ $f\circ \omega^A=\omega^B$ and if for every $n\geq 1$ and $\omega\in \Omega(n)$, the identity:
	
	\begin{equation}
		f(\omega^A(x_1,...,x_n))=\omega^B(f(x_1),...,f(x_n))
	\end{equation}
	holds for every $x_1,...,x_n\in A$.
	
	A vector subspace $I$ of $A$ is an \emph{ideal} of $A$ if for every $n\geq 1$ and $\omega\in \Omega(n)$, we have that:
	
	\[\omega(I,A,\dots,A)+\omega(A,I,\dots,A)+\dots +\omega(A,A,\dots, I)\subseteq I.\]
	
	We say that an $\Omega$-algebra is \emph{simple} if it has no nonzero proper ideals.
	
	Given a ring extension  $R\hookrightarrow S$ and an $R$-linear $\Omega$-algebra $(A,\Omega)$, we can define an $S$-linear $\Omega_S$-algebra, which we denote by $(A_S,\Omega_S)$, by extension of scalar, i.e., it is the algebra whose $S$-module is $A_S=A\otimes_R S$ and such that for every $\omega\in \Omega(n)$, there is an operator $\omega_S\in \Omega_S(n)$ defined by:
	
	\[\omega_S(x_1\otimes s_1,...,x_n\otimes s_n)=\omega(x_1,...,x_n)\otimes s_1\cdots s_n\]
	for every $x_1,...,x_n\in A$ and $s_1,...,s_n\in S$. The definition of monomorphism, empimorphism, automorphism, etc. are as usual. If we have an $R$-linear $\Omega$-algebra $A$, we denote its automorphism group by:
	
	\[\Aut_R(A,\Omega).\]
	
	\begin{remark}\label{rm:MorphismsOmegaDecomp}
		If we have two $R$-linear $\Omega$-algebras $A$ and $B$ and two vector space decompositions  as in Example \ref{ex:ExamplesOmega} (4) indexed by the same set $S$ then, a morphism $\varphi$ between $(A,\Omega^S)$ and $(B,\Omega^S)$ is a morphism between $(A,\Omega)$ and $(B,\Omega)$ satisfying:
		
		\[\varphi(A_s)\subseteq B_s\]
		for all $s\in S$.
	\end{remark}
	
	Given an $\FF$-linear $\Omega$-algebra $(A,\Omega)$, we define its automorphism group scheme as the subgroup scheme of $\GLs(A)$, whose $R$-points are:
	
	\[\AAut(A,\Omega)(R)=\Aut_R(A_R,\Omega_R).\]
	Thus, for a morphism of associative, commutative, unitary algebras $f\colon R\to S$,
	if $\varphi\in \AAut(A,\Omega)(R)$, then, $\AAut(A,\Omega)(f)$ is the $S$-linear extension of the morphism sending $x\otimes 1$ to $(\id\otimes f)\circ \varphi(x\otimes 1)$.

	\subsection{Gradings}
	
	Now, we are going to state some results which appear in the monograph \cite{EKmon} in terms of $\FF$-linear $\Omega$-algebras. The proofs are the same as in the context of algebras so we are not going to write them. For the basic concept about gradings, the reader can refer to the monograph.
	
	\begin{df}
		Let $G$ be a group. A $G$-\emph{grading} $\Gamma$ on an $\FF$-linear $\Omega$-algebra $A$ is a vector space decomposition:
		\[\Gamma\colon A=\bigoplus_{g\in G}A_g\]
		satisfying that for every $\omega\in \Omega(0)$, $\omega(1)=A_e$ and for every $\omega\in \Omega(n)$ with $n\geq 1$,
		\[\omega(A_{g_1},...,A_{g_n})\subseteq A_{g_1\cdots g_n}\]
		for every $g_1,...,g_n\in G$. We say that $(A,\Omega,\Gamma)$ is an $\FF$-linear $G$-\emph{graded} $\Omega$-\emph{algebra}. We call the subspaces $A_g$ the \emph{homogeneous componens of the grading.}
		
		If we have an $H$-grading $\Gamma'$ on $A$:
		\[\Gamma\colon A=\bigoplus_{h\in H}A_h,\]
		we say that $\Gamma'$ is a \emph{coarsening} of $\Gamma$ (or that $\Gamma$ is a \emph{refinement} of $\Gamma'$) if every homogeneous component of $\Gamma$ is contained in a homogeneous component of $\Gamma'$. 
	\end{df}

	\begin{examples}
		\begin{itemize}
			\item[(1)] If $V$ is a left $A$-module for a $G$-graded associative algebra $A$, any decomposition of $V$ as a direct sum of submodules indexed by $G$:
			\[\Gamma \colon V=\bigoplus_{g\in G}V_g\]
			
			satisfying that $a_hV_g\subseteq V_{hg}$ is a $G$-grading of $V$. We say that this is a $G$-graded left $A$-module and concretely, it is a grading of the underlying vector space $V$ (similarly we can define $G$-graded left $A$-modules).
			\item[(2)] In the previous conditions, given $g\in G$, we can define $V^{[g]}$ as the $G$-graded right $A$-module $V$ with grading:
			
			\[\Gamma^{[g]}\colon V=\bigoplus \bar{V}_h\]
			
			where $\bar{V}_h=V_{hg^{-1}}$
			\item[(3)] For a $G$-graded $A$-module $V$ we can consider the associative algebra $\End_{A}(V)$ of those linear morphisms $\varphi\colon V\to V$ satisfying that $\varphi(av)=a\varphi(v)$. We can induce $\End_{A}(V)$ with a grading given by:
			
			\[\End_{A}(V)_g=\{\varphi\mid \varphi(V_h)\subseteq V_{gh}\}\]
		\end{itemize}
	\end{examples}
	
	We say that an $\Omega$-algebra with a $G$-grading $\Gamma$ is graded simple if the corresponding $\Omega^G$-algebra is simple.
	
	We can define a morphism between $\FF$-linear $G$-graded $\Omega$-algebras as a morphism between the corresponding $\Omega^G$ algebras given in Example \ref{ex:ExamplesOmega}. Similarly, we can define monomorphisms, epimorphisms, isomorphisms, automorphisms...
	\begin{remark}
		Notice that due to Remark \ref{rm:MorphismsOmegaDecomp}, $\AAut(A,\Omega,\Gamma)=\AAut(A,\Omega^G)$
	\end{remark}
	If we have a grading as in the previous definition, a group $H$ and a group homomorphism $\alpha\colon G\to H$, we denote by $^{\alpha}\Gamma$ the $H$-grading given by:
	
	\[^{\alpha}\Gamma\colon A=\bigoplus_{h\in H}A_h\]
	where $A_h=\bigoplus_{\underset{g\in G}{\alpha(g)=h}}A_g$.
	
	Given an abelian group $G$, a $G$-grading $\Gamma$ on an $\FF$-linear $\Omega$-algebra $(A,\Omega)$ determines a morphism of affine group schemes:
	
	\[\varphi_{\Gamma}\colon G^D\to \AAut(A,\Omega)\]
	where $G^D$ is the diagonalizable group scheme represented by the group algebra $\FF G$. The morphism is the one induced by the comodule morphism
	
	\[\rho_{\Gamma}\colon A\to A\otimes_\FF \FF G\]
	given by $\rho_{\Gamma}(a)=a\otimes g$ for every $a\in A$ of degree $g$ by the Yoneda's lemma, i.e, $\rho_{\Gamma}=\varphi(\id_{\FF G})$. Conversely, any morphism between a diagonalizable group scheme and the automorphism group scheme of an $\Omega$-algebra induces a grading. An important result for us is the following:
	
	\begin{proposition}\label{p:TransferThm}
		Let $\Omega$ and $\Omega'$ be two operator domains. Let $A$ be an $\Omega$-algebra and $B$ and $\Omega'$ algebra. Assume that we have a morphism:
		\[\theta\colon \AAut(A,\Omega)\to \AAut(B,\Omega')\]
		Then we have a mapping $\Gamma\to \theta(\Gamma)$ from $G$-gradings on $(A,\Omega)$ to $G$-gradings on $(B,\Omega')$ induced by $\varphi_{\Gamma}\mapsto \theta\circ \varphi_{\Gamma}$. If $\Gamma$ and $\Gamma'$ are isomorphic, then $\theta(\Gamma)$ and $\theta(\Gamma')$ are isomorphic.
	\end{proposition}
	\begin{proof}
		The proof is analogous to the proof of \cite[Theorem 1.38]{EKmon}
	\end{proof}
	
	\begin{remark}
		With  the notation of Proposition \ref{p:TransferThm}, if $\theta$ is an isomorphism, then there is a correspondence between the gradings in $(A,\Omega)$ and $(B,\Omega')$ up to isomorphism.
	\end{remark}

	\section{Associative triple systems}\label{se:AT2}
	
	\subsection{Associative envelope}
	
	Our main object of study in this work are associative triple systems of the second kind. What we are going to do now is to give the construction of an associative algebra with an involution and a $3$-grading which was introduced by Ottmar Loos \cite{loos}.

	\begin{df}
		An \emph{associative triple system of the second kind} (AT2) is a ternary algebra $(W,\tripleprod)$ which satisfies that for all $u,v,x,y,z\in W$:
		
		\begin{equation}\label{eq:AT2}
			\{\{u,v,x\},y,z\}=\{u,\{y,x,v\},w\}=\{u,v,\{x,y,z\}\}
		\end{equation}
		
		A subspace $I\subseteq W$ is an \emph{ideal} of $W$ if:

		\begin{equation}\label{eq:AT2}
			\{I,W,W\}+\{W,I,W\}+\{W,W,I\}\subseteq I
		\end{equation}
		
		It is said that $(W,\tripleprod)$ is \emph{simple} if $\{W,W,W\}\neq (0)$ and it has no proper ideals.
	\end{df}
	
	\begin{remark}
		Sometimes, by abuse of notation we will write just $W$ instead of $(W,\tripleprod)$
	\end{remark}
	
	\begin{example}\label{ex:AT2}
		Let $(\cA,\varphi)$ be an associative algebra with involution. We define a triple product on $\cA$ by $\{a,b,c\}=a\varphi(b)c$ for all $a,b,c\in\cA$. This is an AT2 (see \cite{loos}). Moreover, any subspace $W\subseteq \cA$ closed under the triple product is another AT2.
	\end{example}
	
	In \cite{loos} it is shown how all the AT2 arise as in example \ref{ex:AT2}. We are going to show how to construct such associative algebra with involution. Let $(W,\tripleprod)$ be an AT2. Denote $E=\End_{\FF}(W)$. For all $x,y\in W$, we define $l(x,y),r(x,y)\in E$ as:
	
	\begin{equation*}
		l(x,y)z=\{x,y,z\}=r(z,y)x
	\end{equation*}
	
	Let $E^{op}$ be the opposite algebra of $E$ and let $(f,g)\mapsto \overline{(f,g)}\bydef (g,f)$ be the cannonical involution of $E\oplus E^{op}$ and $E^{op}\oplus E$. For all $x,y\in W$ we define:
	
	\begin{align*}
		&\lambda(x,y)\bydef(l(x,y),l(y,x))\in E\oplus E^{op}\\
		&\rho(x,y)\bydef (r(y,x),r(x,y))\in E^{op}\oplus E
	\end{align*}
	
	We denote $L_0$ (resp $R_0$) the subspace of $E\oplus E^{op}$ (resp. $E^{op}\oplus E$) generated by the elements $\lambda(x,y)$ (resp. $\rho(x,y)$) where $x,y\in W$. Finally, denote $e_1$ (resp. $e_2$) the unit of $E\oplus E^{op}$ (resp. $E^{op}\oplus E$).
	
	\begin{lemma}\label{l:AT2Subalgebras}
		$L=\FF e_1\oplus L_0$ (resp. $R=\FF e_2\oplus R_0$) is a subalgebra of $E\oplus E^{op}$ (resp. $E^{op}\oplus E$) which is stable under the involution. Moreover $L_0$ (resp. $R_0$) is an ideal of $(L,\invol)$ (resp. $(R,\invol)$).
	\end{lemma}
	\begin{proof}
		\cite[Lemma 1]{loos}
	\end{proof}
	
	Let $a=(a_1,a_2)\in L$ and $b=(b_1,b_2)\in R$. For $x\in W$, we define:
	
	\begin{equation*}
		ax=a_1(x) \ \ xa=\overline{a}x=a_2(x) \ \ xb=b_1(x) \ \ bx=\overline{b}x=b_2(x)
	\end{equation*}
	
	Let $\overline{W}$ be a copy of $W$. We define $\cA(W)= L\oplus W\oplus \overline{W}\oplus R$.  We can write the elements of $\cA$ in matrix form as:
	
	\begin{equation*}
		\begin{pmatrix}
			a&x\\ y& b
		\end{pmatrix} \ \forall a\in L, x\in W, y\in \overline{W}, b\in R
	\end{equation*}
	
	\begin{theorem}\label{t:AT2ToAlgebras}
		With the previous notation:
		
		\begin{itemize}
			\item[(a)] With the product:
			\begin{equation*}
				\begin{pmatrix}
					a&x\\ y& b
				\end{pmatrix}\begin{pmatrix}
					a'&x'\\ y'& b'
				\end{pmatrix}=\begin{pmatrix}aa'+\lambda(x,y')& ax'+ xb'\\ ya'+by'& \rho(y,x')+bb'\end{pmatrix}
			\end{equation*}
			$\cA(W)$ is an associative algebra with unit $1=\begin{pmatrix}
				e_1&0\\0&e_2
			\end{pmatrix}$ and $\cA(W)_0=L_0\oplus W\oplus \overline{W}\oplus R_0$ is an ideal.
			\item[(b)] $\varphi\colon u=\begin{pmatrix}
				a&x\\ y& b
			\end{pmatrix}\mapsto \overline{u}\bydef \begin{pmatrix}
				\overline{a}& y\\ x& \overline{b}
			\end{pmatrix}$ is an involution of $\cA(W)$ with the product as it is defined in $(a)$ and $\{w_1,w_2,w_3\}=w_1\varphi(w_2)w_3$  $\forall w_1,w_2,w_3\in W$.
			
			\item[(c)] $e_1$ and $e_2$ are orthogonal idempotents of $\cA(W)$ with $\varphi(e_1)=e_1$ and $\varphi(e_2)=e_2$ and they induce a Peirce decomposition:
			
			\begin{equation*}
				\cA(W)_{1,1}=L, \ \cA(W)_{1,2}=W, \ \cA(W)_{2,1}=\overline{W}, \ \cA(W)_{2,2}=R
			\end{equation*}
			
			\item[(d)] $L\oplus R$ doesn't contain any nontrivial ideal of $\cA(W)$
		\end{itemize}
	\end{theorem}
	
	\begin{proof}
		\cite[Theorem 1]{loos}
	\end{proof}
	
	\begin{remark}
		The previous result works for not necessarily finite dimensional AT2.
	\end{remark}
	
	We denote $(\cA(W),\varphi)$ the associative algebra with involution defined in theorem \ref{t:AT2ToAlgebras} and we denote $(\cA(W),\varphi, \Gamma(W))$ the associative algebra with involution $(\cA(W),\varphi)$ fixing the $3$-grading $\Gamma(W)\colon \cA(W)=\cA(W)_{-1}\oplus \cA(W)_0\oplus \cA(W)_1$ where relative to the Peirce decomposition given in theorem \ref{t:AT2ToAlgebras} $\cA(W)_{-1}=\cA(W)_{1,2}$, $\cA(W)_0=\cA(W)_{1,1}\oplus \cA(W)_{2,2}$ and $\cA(W)_1=\cA(W)_{2,1}$. Notice that $\varphi(\cA(W)_1)=\cA(W)_{-1}$.
	
	\begin{remark}
		If $(\cA,\varphi)$ is an algebra with involution and $\Gamma\colon \cA=\cA_{-1}\oplus \cA_0\oplus \cA_1$ is a $3$-grading such that $\varphi(\cA_{i})=\cA_{-i}$ for $i=1,0,-1$, we can define an AT2 $(\cA_{-1},\tripleprod)$ where $\{x,y,z\}\bydef x\varphi(y)z$. We denote this AT2 by $\cW(\cA,\varphi,\Gamma)$.
	\end{remark}
	
	\begin{example}\label{ex:ATtwo}
		In case $\cA=\cM_{n+m}(\FF)$, the $3$-grading $\Gamma$ is given by:
		\begin{align*}&\cA_0=\left\{\begin{pmatrix}
				X&0\\0&Y\end{pmatrix}\mid X\in\cM_n(\FF), Y\in \cM_m(\FF)
			\right\},\\
			&\cA_{-1}=\left\{\begin{pmatrix}
				0&X\\0&0\end{pmatrix}\mid X\in\cM_{n\times m}(\FF)\right\},\\
			&\cA_{-1}=\left\{\begin{pmatrix}
				0&0\\X&0\end{pmatrix}\mid X\in\cM_{n\times m}(\FF)\right\},
		\end{align*}
		and there are matrices $\Phi_1\in \cM_n(\FF)$ and $\Phi_2\in \cM_m(\FF)$ such that:
		\[\Phi=\begin{pmatrix}\Phi_1&0\\0&\Phi_2\end{pmatrix}\]
		satisfies that $\varphi(X)=\Phi X^t\Phi^{-1}$
		then $\cW(\cA,\varphi,\Gamma)$ is isomorphic to $(\cM_{n\times m}(\FF),\tripleprod)$ where $\{X,Y,Z\}=X\Phi_1Y^t\Phi_2^{-1}Z$.
	\end{example}
	
	\begin{example}\label{ex:ATtwoex}
		In case $\cA=\cM_{n+m}(\FF)\oplus \cM_{n+m}(\FF)^{\mathrm{op}}$, the $3$-grading $\Gamma$ is given by:
		\begin{align*}&\cA_0=\left\{\left(\begin{pmatrix}
				X_1&0\\0&Y_1\end{pmatrix},\begin{pmatrix}
				X_2&0\\0&Y_2\end{pmatrix}\right)\mid X_i\in\cM_n(\FF), Y_i\in \cM_m(\FF)
			\right\},\\
			&\cA_{-1}=\left\{\left(\begin{pmatrix}
				0&X_1\\0&0\end{pmatrix},\begin{pmatrix}
				0&X_2\\0&0\end{pmatrix}\right)\mid X_i\in\cM_{n\times m}(\FF)\right\},\\
			&\cA_{-1}=\left\{\left(\begin{pmatrix}
				0&0\\X_1&0\end{pmatrix},\begin{pmatrix}
				0&0\\X_2&0\end{pmatrix}\right)\mid X_i\in\cM_{n\times m}(\FF)\right\},
		\end{align*}
		and  $\varphi(A,B)=(B,A)$
		then $\cW(\cA,\varphi,\Gamma)$ is isomorphic to $(\cM_{n\times m}(\FF)\oplus \cM_{n\times m}(\FF),\tripleprod)$ where $\{(X_1,X_2),(Y_1,Y_2),(Z_1,Z_2)\}=(X_1Y_2Z_1,Z_2Y_1X_2)$.
	\end{example}
	
	\begin{remark}\label{r:AbuseNotationAlgebra}
		Sometimes if the context is clear we will write $\cA$ instead of $(\cA,\varphi,\Gamma)$
	\end{remark}
	
	\begin{remark}\label{r:RecoverAT2}
		Due to theorem \ref{t:AT2ToAlgebras}, we have that $\cW(\cA(W),\varphi,\Gamma(W))=W$
	\end{remark}
	
	Let's see in which conditions the converse is true.
	
	\begin{theorem}\label{t:SimpleSiiSimple}
		An AT2 $W$ is simple if and only if $(\cA(W),\varphi)$ is simple.
	\end{theorem}
	
	\begin{proof}
		\cite[Theorem 3]{loos}
	\end{proof}

	\begin{lemma}\label{l:ZgradingPierce}
		If $(\cA,\varphi)$ is a simple associative algebra (not necessarily finite dimensional) with involution and a non trivial $3$-grading $\Gamma\colon \cA=\cA_{-1}\oplus\cA_0\oplus \cA_1$ such that $\varphi(\cA_i)=\cA_{-i}$ for $i=-1,0,1$, then either $\cA_1=0=\cA_{-1}$ or $\cA_0=\cA_1 \cA_{-1}\oplus \cA_{-1}\cA_1$.
	\end{lemma}
	
	\begin{proof}
		
		Assume that $\cA_1\neq 0$, which after applying the involution is equivalent to $\cA_{-1}$. Consider the subspace \[I=\cA_{-1}\oplus (\cA_1 \cA_{-1}+ \cA_{-1}\cA_1)\oplus \cA_{1}.\] This an ideal of $(\cA,\varphi)$ since for any homogeneous element $a\in\cA$, $aI+Ia\subseteq I$, and clearly $\varphi(I)=I$. Since $I$ is a nonzero ideal and this algebra with involution is simple, it follows that $\cA=I$. Hence $\cA_0=\cA_1 \cA_{-1}+ \cA_{-1}\cA_1$.
		
		Let \[J=(\cA_1 \cA_{-1})\cap (\cA_{-1}\cA_1).\] Due to the fact that \[(A_1+A_{-1})J+J(A_1+A_{-1})=(0)\] and \[A_0J+JA_0\subseteq J,\] $J$ is another ideal. Since it is contained in $\cA_0$, it is not $\cA$. Thus, $J=(0)$. Therefore, the sum is direct, i.e, $\cA_0=\cA_1 \cA_{-1}\oplus \cA_{-1}\cA_1$.
	\end{proof}
	
	\begin{lemma}\label{l:ASimpleWSimple}
		If $(\cA,\varphi)$ is a simple associative algebra with involution and a $3$-grading $\Gamma\colon \cA=\cA_{-1}\oplus\cA_0\oplus \cA_1$ such that $\varphi(\cA_i)=\cA_{-i}$ for all $i=-1,0,1$, then either $\cW(\cA,\varphi,\Gamma)$ is simple or 0.
	\end{lemma}
	
	\begin{proof}
		Call $W=\cW(\cA,\varphi,\Gamma)$. We will prove by contradiction that $\{W,W,W\}\not=(0)$. Arguing by contradiction, if $\{W,W,W\}=(0)$, then by definition $\cA_{-1}\varphi(\cA_{-1})\cA_{-1}=\cA_{-1}\cA_1\cA_{-1}=(0)$. This fact and the fact that $\cA_1^2=\cA_{-1}^{2}=(0)$ would mean that $\cA_0=\cA_1\cA_{-1}\oplus \cA_{-1}\cA_1$ is an ideal of $\cA$. Since $\varphi(\cA_0)=\cA_0$, it follows $\cA_0$ is an ideal of $(\cA,\varphi)$ and by simplicity $\cA_0=(0)$. Hence $\cA^2=(0)$ and this is a contradiction with the fact that $(\cA,\varphi)$ is simple.
		
		Let $I$ be a nonzero ideal of $W$, we claim that \[J=I\oplus(I\cA_1+\cA_1 I+\varphi(I)\cA_{-1}+\cA_{-1}\varphi(I))\oplus\varphi(I)\]
		is an ideal of $(\cA,\varphi)$. It is clearly stable under the involution. Moreover by definition of ideal in an AT2, the inclusion \[I\cA_{1}\cA_{-1}+\cA_{-1} \varphi(I)\cA_{-1}+\cA_{-1}\cA_1I\subseteq I\] holds. Thus, it is easy to prove that $\cA_i J+J\cA_i\subseteq J$  for $i=-1,1$ and since due to Lemma \ref{l:ZgradingPierce}, $\cA_0=\cA_1 \cA_{-1}\oplus \cA_{-1}\cA_1$, it follows that $\cA_0 J+J\cA_0\subseteq J$ implying that  $J$ is an ideal. Hence $J=\cA$ and $J\cap \cA_{-1}=I$. Therefore, $I=W$.
	\end{proof}
	
	\begin{lemma}\label{l:AssocAlgIso}
		If $(\cA,\varphi)$ is a simple associative algebra with involution and with a $3$-grading $\Gamma\colon \cA=\cA_{-1}\oplus \cA_0\oplus \cA_1$ such that $\varphi(\cA_i)=\cA_{-i}$ and $\cA_{-1}\neq 0$, then \[(\cA,\varphi,\Gamma)\cong(\cA(\cW(\cA)),\varphi,\Gamma(\cW(\cA)))\]
	\end{lemma}
	
	\begin{proof}
		We define $\psi\colon \cA\to \cA(\cW(\cA)))$ as the vector space homomorphism induced by $\psi(a)=a$ for all $a\in \cA_{1}$, $\psi(b)=\varphi(\psi(\varphi(b)))$ for all $b\in \cA_1$ and $\psi(x)=\sum_{i=1}^k\psi(a_i)\psi(b_i)$ if $x=\sum_{i=1}^ka_ib_i$ with all $a_i\in \cA_{-1}$ and all $b_i\in \cA_{1}$ or viceversa. This is enough due to lemma \ref{l:ZgradingPierce}.
		
		In order to prove that it is well defined we need to prove that if $x=\sum_{i=1}^k a_ib_i$ with all $a_i\in \cA_{-1}$ and all $b_i\in \cA_{1}$ or viceversa, then $\sum_{i=1}^k\psi(a_i)\psi(b_i)=0$. Suppose that all $a_i\in \cA_{-1}$ and all $b_i\in \cA_{1}$. Because of theorem \ref{t:AT2ToAlgebras} $\sum_{i=1}^k\psi(a_i)\psi(b_i)=(x_1,x_2)\in \End_{\FF}(\psi(\cA_{-1}))\oplus \End_{\FF}(\psi(\cA_{-1}))^{op}$. Let $a\in \cA_{-1}$:
		
		\[x_1(a)=\psi(x)a=\sum_{i=1}^k\psi(a_i)\psi(b_i)a=\psi(a_i)\varphi(\psi(\varphi(b_i)))\psi(a)=\sum_{i=1}^ka_ib_ia=0\]
		
		Hence $x_1=0$. Similarly we do with $x_2$ taking $a\in\cA_1$ and using that $a(x_1,x_2)=(x_2,x_1)a$. When  all $a_i\in \cA_{1}$ and all $b_i\in \cA_{-1}$ we do it in the same way.
		
		In order to prove that it is an algebra homomorphism we notice that we only have to check that $\psi(a)\psi(b)=\psi(ab)$ for homogeneous elements. If $a\in \cA_{-1}$ and $b\in \cA_1$, or viceversa, the identity is true by definition. If $a$ or $b$ are in $\cA_0$, using lemma \ref{l:ZgradingPierce} we can suppose that they are in $\cA_1\cA_{-1}$ or in $\cA_{-1}\cA_1$, hence, we just need to prove that $\psi(a)\psi(b)\psi(c)=\psi(abc)$ and $\psi(a)\psi(b)\psi(c)\psi(d)=\psi(abcd)$ with $a,b,c,d\in \cA_1\cup \cA_{-1}$ and this is a direct calculation.
		
		$\psi$ preserves the involution, because it preserves the involution in $\cA_1\oplus \cA_{-1}$ and due to lemma \ref{l:ZgradingPierce}.
		
		Since both algebras are generated as algebras with involution by he homogeneous components of degree $-1$ (lemma \ref{l:ZgradingPierce}) $\psi$ is surjective.
		
		Since $\ker\psi$ is an ideal of $(\cA,\varphi)$, $(\cA,\varphi)$ is simple and $\ker \psi \cap (\cA_1\oplus \cA_{-1})=(0)$, we get that $\ker \psi=(0)$.
	\end{proof}
	
	We will finish with a lemma.

	\begin{lemma}\label{l:PierceCS}
		Let $(\cA,\varphi)$ be a central simple associative algebras with involution over a field $\FF$. Let $e$ be an idempotent element such that $\varphi(e)=e$ and let $\cA_{1,1}=e\cA e$. Then $\varphi(\cA_{1,1})=\cA_{1,1}$ and $(\cA_{1,1},\varphi_{\mid\cA_{1,1}})$ is central simple
	\end{lemma}
	\begin{proof}
		$\varphi(\cA_{1,1})=\cA_{1,1}$ is clear. 
		
		For the second part, without loss of generality, we can suppose that $\FF$ is algebraically closed. 
		
		If $\cA$ is central simple, then, there is a vector space $V$ such that  $A\cong \End_{\FF}(V)$. Let $e\in\End_{\FF}(V)$ be an idempotent. Then, $\mathrm{Im}(e)=\{v\in V\mid e(v)=v\}$ and $V=\mathrm{Im}(e)\oplus\ker(e)$.
		The morphism $\psi \colon e\End_{\FF}(V) e\to \End_{\FF} \mathrm{Im}(e)$ given by $\psi(f)=f_{\mid \mathrm{Im}(e)}$ is an isomorphism, Therefore, $e\End_{\FF}(V) e=\End_{\FF}(\mathrm{Im}(e))$ which is simple. Therefore we get the result for this case.
		
		If $\cA$ is not central simple, then, there is a central simple algebra $\cB$ with $(\cA,\varphi)\cong(\cB\oplus \cB^{op},\ex)$ where $\ex(b_1,b_2)=(b_2,b_1)$. If $e\in \cB\oplus \cB^{op}$, there is $e_1,e_2\in\cB$ such that $e=(e_1,e_2)$. Since $\ex(e)=e$, we get that $e_1=e_2$. Since $e^2=e$ we get that $e_1^2=e_1$. Therefore, $e(\cB\oplus \cB^{op}) e=e_1\cB e_1\oplus (e_1\cB e_1)^{op}$. Since $e_1\cB e_1$ is central simple and $\ex$ restricts to $e(\cB\oplus \cB^{op}) e$, we have that $(\cA_{1,1},\varphi_{\mid\cA_{1,1}})$  is central simple.
	\end{proof}
	\subsection{Automorphism group schemes}
	
	On this section we are going to find an isomorphism between the automorphism group scheme of an AT2 $W$ and the automorphism group  scheme of the corresponding associative algebra with an involution and a $3$-grading $(\cA(W),\varphi,\Gamma)$. We will assume that the reader has a basic knowledge about affine group scheme but in order to get an introduction to the topic, the following references are useful \cite{EKmon}, \cite{KMRT} and \cite{WatGroupSch}.
	
	For an AT2 $(W,\tripleprod)$, we define its automorphism group scheme $\AAut(W,\tripleprod)$ as the affine group scheme whose group of $R$-points is:
	
	\[\AAut(W,\tripleprod)(R)=\Aut_R(W_R,\tripleprod_R).\]
	
	For an associative algebra $\cA$, with an involution $\varphi$ and a $3$-grading $\Gamma\colon \cA=\cA_{-1}\oplus\cA_0\oplus \cA_1$, we define its automorphism group scheme $\AAut(\cA,\varphi,\Gamma)$ as the affine group scheme whose group of $R$-points is:
	
	\[\AAut(\cA,\varphi,\Gamma)(R)= \Aut_R(\cA_R,\varphi_R,\Gamma_R).\]

	Our purpose will be to prove that for any simple AT2 $W$,  there is an isomorphism between $\AAut(W,\tripleprod)$ and $\AAut(\cA(W),\varphi,\Gamma(W))$. In order to so, we should start by proving a technical lemma.
	
	\begin{lemma}\label{l:AutomorphismWellDeff}
		Let $(W,\tripleprod)$ be an AT2, let $R$ be an unital, commutative, associative algebra and let $\psi\in \AAut(W,\tripleprod)(R)$. Given $x_1,...,x_k, y_1,...,y_k\in W_R$, $\sum_{i=1}^{k}l(x_i,y_i)=0$ (resp. $\sum_{i=1}^{k}r(x_i,y_i)=0$) if and only if $\sum_{i=1}^{k}l(\psi(x_i),\psi(y_i))=0$  (resp. $\sum_{i=1}^{k}r(\psi(x_i),\psi(y_i))=0$).
	\end{lemma} 
	\begin{proof}
		We are going to prove the lemma for $l$. The proof for $r$ it is analogous. Since $\psi$ is surjective, given $w\in W$, there is $v\in W$ such that $\psi(v)=w$. Now since $\psi$ is an automorphism of the AT2:
		
		\[\sum_{i=1}^{k}l(\psi(x_i),\psi(y_i))w=\sum_{i=1}^{k}l(\psi(x_i),\psi(y_i))\psi(v)=\psi(\sum_{i=1}^{k}l(x_i,y_i)v)=0.\]	 
		
		The fact that $\psi$ is injective implies the statement.
	\end{proof}
	
	\begin{remark}\label{rm:AutomorphismWellDeff}
		By definition of $\lambda$ (resp. $\rho$), due to Lemma \ref{l:AutomorphismWellDeff}, for an automorphism $\psi\in \AAut(W,\tripleprod)(R)$, if $x_1,...,x_k, y_1,...,y_k\in W_R$, it follows that $\sum_{i=1}^{k}\lambda(x_i,y_i)=0$ (resp. $\sum_{i=1}^{k}\rho(x_i,y_i)=0$) implies that $\sum_{i=1}^{k}\lambda(\psi(x_i),\psi(y_i))=0$  (resp. $\sum_{i=1}^{k}\rho(\psi(x_i),\psi(y_i))=0$).
	\end{remark}
	
	\begin{remark}\label{rm:HomomorphismDef}
		Let $(W,\tripleprod)$ be a simple AT2, let $R$ be an unital, commutative, associative algebra and let $\psi\in \AAut(W,\tripleprod)(R)$. We can define a vector space homomorphism:
		
		\[\cA(\psi)\colon \cA(W_R)\to \cA(W_R)\]
		by $\cA(\psi)(a)=a$ if $a\in\cA(W_R)_{-1}$, $\cA(\psi)(a)=\varphi(\psi(\varphi(a)))$ if $a\in\cA(W_R)_1$ and $\cA(\psi)(a)=\sum_{i=1}^k\psi(a_i)\varphi(\psi(\varphi (b_i)))$ if $a=\sum_{i=1}^ka_ib_i$ with $a_i\in\cA(W_R)_1$ and $b_i\in\cA(W_R)_{-1}$ or viceversa. This is well defined due to the fact that since $\cA(W)$ is simple, Lemma \ref{l:ZgradingPierce} implies that \[\cA(W_R)_0=\cA(W_R)_1\cA(W_R)_{-1}\oplus\cA(W_R)_{-1}\cA(W_R)_1\] and due to 
		Remark \ref{rm:AutomorphismWellDeff}, if $a=0$, $\psi(a)=0$.
		
		Since $\cA(\psi)$ preserves the degrees it follows that it preserves the involution. 
		
		In order to show that it preserves the product, we can just take homogeneous elements. Concretely, for elements from $\cA(W_R)_0$, we can just choose elements of the form $ab$ with $a\in \cA(W_R)_1$ and $b\in \cA(W_R)_{-1}$ or viceversa. Then:
		
		\begin{itemize}
			\item	If $a\in \cA(W_R)_1$ and $b\in \cA(W_R)_{-1}$, it follows by definition that $\cA(\psi)(ab)=\psi(a)\varphi(\psi(\varphi(b)))=\cA(\psi)(a)\cA(\psi(b))$.
			\item  If $a\in \cA(W_R)_{-1}$ and $b\in \cA(W_R)_{1}$ we do as before.
			\item  If $a,b\in \cA(W_R)_{1}$ and $c\in \cA(W_R)_{-1}$ it follows that $ab=0$, and so, 
			\begin{equation*}
				\begin{split}
					\cA(\psi)(a(bc))=\cA(\psi)((ab)c)=0=(\cA(\psi)(a)\cA(\psi)(b))\cA(\psi)(c)\\=\cA(\psi)(a)(\cA(\psi)(b)\cA(\psi)(c))=\cA(\psi(a)))\cA(\psi(bc)))
				\end{split}
			\end{equation*}
			
			And on the other hand:
			
			\begin{equation*}
				\begin{split}	\cA(\psi)(a(cb))=\cA(\psi)(\{a,\varphi(c),b\})=\{\psi(a),\psi(\varphi(c)),\psi(b)\}\\=\psi(a)\varphi(\psi(\varphi(c))\psi(b)=\cA(\psi)(a)\cA(\psi)(c)\cA(\psi)(b)
				\end{split}
			\end{equation*}
			\item   If $a,b\in \cA(W_R)_{-1}$ and $c\in \cA(W_R)_{1}$ the equatities  	$\cA(\psi)(a(bc)) =\cA(\psi(a)))\cA(\psi(bc)))$ and $\cA(\psi)(a(cb))=\cA(\psi)(a)\cA(\psi)(c)\cA(\psi)(b)$ follow as in the previous item.
			
			\item   If $a,c\in \cA(W_R)_{1}$ and $b,d\in \cA(W_R)_{-1}$ we have that:
			\begin{equation*}
				\begin{split}
					\cA(\psi)((ab)(cd))=\cA(\psi)(((ab)c)d)=\cA(\psi)((ab)c))\cA(\psi)(d)\\=\cA(\psi)(ab)\cA(\psi)(c)(\psi)(d)=\cA(\psi)(ab)\cA(\psi)(cd)
				\end{split}
			\end{equation*}
			And the fact that it preserves the other possible products is proved similarly.
		\end{itemize}

		Finally, we have that if $\psi,\varphi\in \AAut(W,\tripleprod)(R)$, it follows that 
		
		\begin{equation}\label{eq:cAComposition}
			\cA(\psi\circ\varphi)=\cA(\psi)\circ \cA(\varphi).
		\end{equation} 
		
		Therefore, choosing $\psi^{-1}$ as $\varphi$, we check that $\cA(\psi)\in \AAut(\cA(W),\invol,\Gamma(W))(R)$.
	\end{remark}
	
	\begin{theorem}\label{t:AutSchemeIso}
		Let $(W,\tripleprod)$ be a simple AT2. With the notation of Remark \ref{rm:HomomorphismDef}, the morphism:
		
		\[\theta\colon \AAut(W,\tripleprod)\to \AAut(\cA(W),\invol,\Gamma(W))\]
		given by $\theta_R(\psi)=\cA(\psi)$ is an isomorphism with inverse  $\theta^{-1}(\psi)=\psi_{\lvert \cA(W_R)_{-1}}$.
	\end{theorem}
	\begin{proof}
		Due to Theorem \ref{t:SimpleSiiSimple} $(\cA(W),\invol)$ is simple. Due to Lemma \ref{l:ZgradingPierce}, $\cA(W)$ is generated as an algebra with involution by $W$. Thus, for every commutative, associative unital algebra $R$, $\cA(W)_R$ is generated as an algebra with involution by $W_R$.
		
		Given $\psi \in \AAut(W,\tripleprod)(R)$, $\cA(\psi)$ is a well defined automorphism due to Remark \ref{rm:HomomorphismDef}. From the definition it follows that $\cA(\psi)_{\lvert \cA(W_R)_{-1}}=\psi$. Now, if $\psi\in \AAut(\cA(W),\invol,\Gamma(W))$, the fact that $\cA(W)$ is generated as an algebra with involution by $W$, implies that $\cA(\psi_{\lvert \cA(W_R)_{-1}})=\psi$. Therefore, the result follows.
	\end{proof}
	
	\section{Gradings on the associative algebras over arbitrary fields of characteristic other than 2}\label{sec:ArbitraryFields}
	
	\subsection{Associative algebra}
	
	Let $G$ be a finitely generated group. We would like to classify up to isomorphism gradings on central simple associative algebras with involution  $(\cA,\varphi)$ with a $3$-grading $\Gamma\colon \cA=\cA_{-1}\oplus \cA_0\oplus \cA_1$ such that $\varphi(\cA_{i})=\cA_{-i}$, i.e., we want to classify up to isomorphism $\Omega$-algebras $(\cA,\varphi,\Gamma_1,\Gamma_2)$ satisfying:
	
	\begin{itemize}
		\item[(Q1)] $(\cA,\varphi)$ is a central simple associative algebra with involution,
		\item[(Q2)] $\Gamma_1$ is a $\ZZ$-grading of $\cA$,
		\item[(Q3)]$\supp(\Gamma_1)= \{-1,0,1\}$ and $\varphi(\cA_i)=\cA_{-i}$ for all $i\in \{-1,0,1\}$
		\item[(Q4)] $\Gamma_2$ is a $G$-grading of $(\cA,\varphi,\Gamma_1)$.
	\end{itemize}  
	
	If we denote the homogeneous components by \[\Gamma_1\colon \cA=\cA_{-1}\oplus \cA_0\oplus \cA_1,\] and \[\Gamma_2\colon \cA=\bigoplus_{g\in G}\cA_g,\]
	we can define a $\ZZ\times G$-grading 
	\[\Gamma\colon \cA=\bigoplus_{(i,g)\in\ZZ\times G}\cA_{(i,g)}\]
	where $\cA_{(i,g)}=\cA_i\cap \cA_g$ for all $(i,g)\in \ZZ\times G$. Denote $\pi_1\colon \ZZ\times G\to \ZZ$ and $\pi_2\colon \ZZ\times G\to G$ the projections on $\ZZ$ and $G$. Then, $(A,\varphi,\Gamma)$ satisfies:
	
	\begin{itemize}
		\item[(T1)] $(A,\varphi)$ is a simple associative algebra with involution,
		\item[(T2)] $\Gamma$ is a $\ZZ\times G$-grading of the algebra $A$,
		\item[(T3)]	$\supp ^{\pi_1}\Gamma= \{-1,0,1\}$, and
		\item[(T4)] $\varphi(A_{(i,g)})=A_{(-i,g)}$ for all $i\in \{1,2,3\}$ and $g\in G$.
	\end{itemize}
	
	Moreover, if we have an $\FF$-lineal $\Omega$-algebra $(A,\varphi,\Gamma)$ satisfying $(T1)-T(4)$, it follows that $(A,\varphi,^{\pi_1}\Gamma,^{\pi_2}\Gamma)$. Additionally, we have that for two pairs of algebras $(A_1,\varphi_1,\Gamma_1)$ and $(A_2,\varphi_2,\Gamma_2)$:
	
	\[\Hom((A_1,\varphi_1,\Gamma_1),(A_2,\varphi_2,\Gamma_2))=\Hom((A_1,\varphi_1,^{\pi_1}\Gamma_1,^{\pi_2}\Gamma_1),(A_2,\varphi_2,^{\pi_1}\Gamma_2,^{\pi_2}\Gamma_2)).\]
	
	Therefore, the problem of classifying up to isomorphism the $\Omega$-algebras $(A,\varphi,\Gamma)$ satisfying $(Q1)-(Q4)$ is equivalent to the problem of classifying up to isomorphism the $\Omega$-algebras $(A,\varphi,\Gamma_1,\Gamma_2)$ satisfying $(T1)-(T4)$. If $(A,\varphi,\Gamma)$ is such kind of algebra, we can reduce the problem into two cases: the case where $(A,\Gamma)$ is graded simple and the case where $(A,\Gamma)$ is not graded simple.
	
	\subsection{Case with $(A,\Gamma)$ not simple}
	
	We will use the following notation: for an algebra $I$, we denote by $\ex$ the \emph{exchange involution} on $A=I\oplus I^{op}$ i.e, 
	\[\ex(x,y)=(y,x).\]
	
	\begin{df}
		Given a $\ZZ\times G$-grading $\Delta$ on $I$, we can define a $\ZZ\times G$-grading on $A=I\oplus I^{op}$, which we call the \emph{exchange grading},  by:
		
		\begin{equation*}\label{eq:exchangeGrading}\Delta^{\ex}\colon A=\bigoplus_{(i,g)\in\ZZ\times G}A_{(i,g)}
		\end{equation*}
		given  by:
		
		\[A_{(i,g)}=\{(x,y)\in A\mid x\in I_{(i,g)}, y\in I_{(-i,g)}\}.\]
		
	\end{df}
	
	\begin{lemma}
		Let $(A,\varphi,\Gamma)$ be an $\Omega$-algebra satisfying $(T1)-(T4)$ such that $(A,\Gamma)$ is not graded simple. Then there is a simple algebra $I$ and a grading $\Gamma_I$ with $\supp(^{\pi_1}\Gamma_I)=\{-1,0,1\}$ such that:
		
		\[(A,\varphi,\Gamma)\cong(I\oplus I^{op},\ex,\Gamma_I^{\ex})\]
	\end{lemma}
	\begin{proof}
		The fact that $(A,\Gamma)$ is not graded simple implies that there is a  proper graded ideal $I\neq 0$ of $(A,\Gamma)$. The ideal $J_1=I+ \varphi(I)$ is a nontrivial ideal of $(A,\varphi)$. Since $(A,\varphi)$ is simple, it implies that $J=A$. The ideal $J_2=I\cap \varphi(I)$ is a proper ideal of $A$. Therefore, $J_2=0$. Thus, $A=I\oplus \varphi(I)$.  Denote 
		\[\Gamma_I\colon I=\bigoplus_{(i,g)\in\ZZ\times G}I_{(i,g)}\]
		the grading induced by the restriction, i.e, 
		\[I_{(i,g)}=I\cap A_{(i,g)}.\]
		The graded morphism 
		
		\[	I\oplus \varphi(I)\to I\oplus I^{op}\]
		given by $x+\varphi(y)\mapsto (x,y)$ for every $x,y\in I$ is an isomorphism. 
		
		Finaly, $I$ is simple because if $0\neq J$ is an ideal of $I$, then, the fact that $A=I\oplus \varphi(I)$ as an algebra imples that $I$ is an ideal of $A$. Then, as shown before, it must imply that $A=J\oplus \varphi(J)$. Thus $I=J$.
	\end{proof}
	
	Due to the previous lemma, it would be enough to classify the graded algebras $(I\oplus I^{op},\ex,\Gamma_I )$ where $I$ is a simple associative algebra and $\Gamma_I$ is a $\ZZ\times G$-grading of $I$ such that $\supp(^{\pi_1}\Gamma_I)=\{-1,0,1\}$.
	
	\begin{df}
		\begin{itemize}
			\item[(1)] Let $G$ be an abelian group and let $D$ be a $G$-graded division algebra with support $T$, let $k$ be a positive integer and $\gamma=(g_1,...,g_k)\in D^k$. We denote $\Gamma(G,D,\gamma)$ the grading on $\cM_n(D)$ given by $\deg(d E_{i,j})=g_i(\deg(d))g_j^{-1}$ for every $d\in D$ and $1\leq i,j\leq k$.
			
			\item[(2)] Let $G$ be a group and $D$ a $\ZZ\times G$-graded division algebra with support $T$ (notice that since $D$ is a graded division algebra of finite dimension, $T\subseteq \{0\}\times G$). Let $g_1,...,g_n\in G$ and let $k_0,k_1$ be two positive integers such that $k_0+k_1=n$. Denote $\gamma_0=(g_{1},...,g_{k_0})$ and $\gamma_1=(g_{k_0+1},...,g_{k_0+k_1})$. We denote by $\Gamma(G,D,\gamma_0,\gamma_1)$ the $\ZZ\times G$-grading on $\cM_n(D)$ given by $\deg(d)E_{i,j}=(\delta_i-\delta_j,g_i\deg(d)g_j^{-1})$ where $\delta_k=0$ if $k\leq k_0$ and $\delta_k=1$ otherwise. We denote:
			\[\cM(G,D,\gamma_0,\gamma_1)=(\cM_n(D),\Gamma(G,D,\gamma_0,\gamma_1))\]
			\item[(3)] With the previous notation, we denote:
			\[\cM(G,D,\gamma_0,\gamma_1)^{\ex}=(\cM_n(D)\oplus \cM_n(D)^{op},\ex,\Gamma(G,D,\gamma_0,\gamma_1)^{\ex}).\]
			\item[(4)] Let $\gamma=(g_1,..,g_n)\in G^{n}$, $\kappa=(\kappa_1,...,\kappa_n)\in\ZZ_{>0}^n$ and let $m=\kappa_1+...+\kappa_n$. We denote $\bold{\gamma}(\kappa,\gamma)=(h_1,...,h_m)\in G^m$ where $h_i=g_j$ if $\kappa_1+...+\kappa_{j-1}<i\leq \kappa_1+...+\kappa_{j}$ for every $i\in\{1,...,m\}$ and $j\in \{1,...,n\}$.
			\item[(5)] Finally, given an $n$-uple $\gamma_0\in G^n$ and a $m$-uple $\gamma_1\in G^m$ consisting on different elements of $G$, we denote:
			
			\begin{align*}
				&\Gamma(G,D,\kappa_0,\gamma_0)=\Gamma(G,D,\bold{\gamma}(\kappa_0,\gamma_0))\\
				&\Gamma(G,D,\kappa_0,\gamma_0,\kappa_1,\gamma_1)=\Gamma(G,D,\bold{\gamma}(\kappa_0,\gamma_0),\bold{\gamma}(\kappa_1,\gamma_1))\\
				&\cM(G,D,\kappa_0,\gamma_0,\kappa_1,\gamma_1)=\cM(G,D,\bold{\gamma}(\kappa_0,\gamma_0),\bold{\gamma}(\kappa_1,\gamma_1))\\
				&\cM(G,D,\kappa_0,\gamma_0,\kappa_1,\gamma_1)^{\ex}=\cM(G,D,\bold{\gamma}(\kappa_0,\gamma_0),\bold{\gamma}(\kappa_1,\gamma_1))^{\ex}
			\end{align*}
		\end{itemize}
	\end{df}
	
	\begin{remark}
		We need to recall that for
		\[\gamma_0=(g_1,...,g_{k_0})\in G^{k_0} \text{ and } \gamma_1=(g_{k_0+1},...,g_{k_0+k_1})\in G^{k_1},\] the $G$-graded algebra $\cM(G,D,\gamma_0,\gamma_1)$ is isomorphic to as a graded algebra with $\End_{D}(V)$ the grading on $V$ is determined by a homogeneous basis $v_1,...,v_{k_0+k_1}$ whose degrees are given by $\deg(v_i)=(\delta_i,g_i)$ with $\delta_i=0$ if $1\leq i\leq k_0$ and $\delta_i=1$ otherwise.
	\end{remark}
	\begin{df}Let $G$ be a group $D$ a $G$-graded division algebra with support $T$. Let $\gamma=(g_1,...,g_n)\in G^n$ be an $n$-uple of elements of $G$. We denote by $\Xi(\gamma)$ the multiset $\{g_1T,...,g_n T\}$ with multiplicity for $g_iT$ the number of appearances of  $g_i$ modulo $T$ in $\gamma$. For $\gamma$ and $n$-uple of distinct elements of $G$ and $\kappa\in \ZZ^n_{>0}$ we denote $\Xi(\kappa,\gamma)=\Xi(\bold{\gamma}(\kappa,\gamma))$.
	\end{df}
	
	\begin{lemma}
		Let $A$ be a simple associative algebra and let $\Gamma$ be a $\ZZ\times G$-grading with $\supp(^{\pi_1}\Gamma)= \{-1,0,1\}$. In this situation, $(A,\Gamma)$ is isomorphic to $\cM(G,D,\gamma_0,\gamma_1)$ for some $\ZZ\times G$-graded division algebra $D$, $\gamma_0=(g_1,...,g_{k_0})\in G^{k_0}$ and $\gamma_1=(g_{k_0+1},...,g_{k_0+k_1})\in G^{k_1}$ and with $k_0,k_1\neq 0$. Moreover, the  $\ZZ\times G$-graded algebras $\cM(G,D,\gamma_0,\gamma_1)$ and $\cM(G,D',\gamma'_0,\gamma'_1)$ are isomorphic if and only if $D$ is isomorphic to $D'$ as graded algebras and there is an element $g\in G$ such that $\Xi(\gamma_i)=g\Xi(\gamma'_i)$ for $i\in\{0,1\}$.
	\end{lemma}
	\begin{proof}
		This proof is based on the result given in \cite{HSK19} in which the case over an algebraically closed field is solved.
		
		As shown in \cite[Theorem 2.10]{EKmon}, $(A,\Gamma)$ is isomorphic to $\End_{D}(V)$ for some graded division algebra $D$ and some graded right $D$-module $V$. Denote $T=\supp(D)$. The fact that $T\subseteq \{-1,0,1\}\times G$ and it is a finite group, implies that $T\subseteq \{0\}\times G$. Due to the fact that $\supp (\Gamma)\subseteq \{-1,0,1\}\times G$, it follows that there is an integer $i$ such that $\supp(V)\subseteq \{i,i+1\}$. \cite[Theorem 2.10]{EKmon} implies that we can do a shift on the grading of $V$, i.e., multiplying every degree by an element $(j,g)\in \ZZ\times G$ and obtain the same grading on $\End_{D}(V)$. Thus, we can assume that $i=0$. Moreover, since $\supp(^{\pi_1}\Gamma)=\{-1,0,1\}$, we can assume that there is one element of $V$ with degree on $\{0\}\times G$ and one element with degree on $\{1\}\times G$.
		
		Taking a homogeneous basis $v_1,...,v_{k_0+k_1}$ of degrees $\deg(v_i)=(\delta_i,g_i)$ with $\delta_i=0$ if $1\leq i\leq k_0$ and $\delta_i=1$ otherwise, we can see that $(A,\Gamma)$ is isomorphic to $\cM(G,D,\gamma_0,\gamma_1)$ where $\gamma_0=(g_1,...,g_{k_0})\in G^{k_0}$ and $\gamma_1=(g_{k_0+1},...,g_{k_0+k_1})\in G^{k_1}$. Finally, since $A$ is simple, $D$ is simple.
		
		The second part follows directly from \cite[Corollary 2.12]{EKmon}.
	\end{proof}
	
	The following result is straightforward from the lemma:
	
	\begin{corollary}
		Let $A$ be a simple associative algebra and let $\Gamma$ be a $\ZZ\times G$-grading with $\supp(^{\pi_1}\Gamma)= \{-1,0,1\}$. In this situation, $(A,\Gamma)$ is isomorphic to $\cM(G,D,\kappa_0,\kappa_1,\gamma_0,\gamma_1)$ for some $\ZZ\times G$-graded division algebra $D$, $\gamma_0=(g_1,...,g_{k_0})\in G^{k_0}$ and $\gamma_1=(g_{k_0+1},...,g_{k_0+k_1})\in G^{k_1}$ consisting on different elements, $\kappa_0\in \ZZ_{>0}^{k_0}$, $\kappa_1\in \ZZ_{>0}^{k_1}$ and with $k_0,k_1\neq 0$. Moreover, the graded algebras $\cM(G,D,\kappa_0,\kappa_1,\gamma_0,\gamma_1)$ and $\cM(G,D',\kappa'_0,\kappa'_1\gamma'_0,\gamma'_1)$ are isomorphic if and only if $D$ is isomorphic to $D'$ as graded algebras and there is an element $g\in G$ such that $\Xi(\kappa_i,\gamma_i)=g\Xi(\kappa'_i,\gamma'_i)$ for $i\in\{0,1\}$.
	\end{corollary}
	
	\begin{df}
		Given an associative algebra $I$, and a $\ZZ\times G$-grading \[\Gamma\colon I=\bigoplus_{(i,g)\in\ZZ\times G}I_{(i,g)},\] we define the $(\ZZ\times G)$-grading:
		\[\Gamma^{op}\colon I^{op}=\bigoplus_{(i,g)\in\ZZ\times G}I^{op}_{(i,g)}\]
		on the opposite algebra $I^{op}$ by $I^{op}_{(i,g)}=I_{(-i,g)}$ for every $(i,g)\in \ZZ\times G$.
	\end{df}
	
	\begin{lemma}
		Let $I$ and $J$ be simple associative algebras. Let $\Gamma_I$ and $\Gamma_J$ be $(\ZZ\times G)$-gradings on $I$ and $J$ respectively satisfying $\supp(\Gamma_I)=\{-1,0,1\}=\supp(\Gamma_J)$. Then, the $\Omega$-algebra $(I\oplus I^{op},\ex,\Gamma_I^{\ex})$ is isomorphic to $(J\oplus J^{op},\ex,\Gamma_J)$ if and only if $(I,\Gamma_I)$ is isomorphic to either $(J,\Gamma_J)$ or to $(J^{op},\Gamma_J^{op})$.
	\end{lemma}
	\begin{proof}
		It is clear that if $(I,\Gamma_I)$ is isomorphic to either $(J,\Gamma_J)$ or to $(J^{op},\Gamma_J^{op})$, then $(I\oplus I^{op},\ex,\Gamma_I^{\ex})$ is isomorphic to $(J\oplus J^{op},\ex,\Gamma_J)$.
		
		Assume that $(I\oplus I^{op},\ex,\Gamma_I^{\ex})$ is isomorphic to $(J\oplus J^{op},\ex,\Gamma_J)$ via an isomorphism
		\[\psi\colon I\oplus I^{op}\to J\oplus J^{op}.\]
		
		In this situation, $\psi(I)$ is a proper nontrivial ideal of $J\oplus J^{op}$. Therefore, $\psi(I)$ is either $J$ or $J^{op}$. In the first case, $(I,\Gamma_I)$ is isomorphic to either $(J,\Gamma_J)$ and in the second case $(I,\Gamma_I)$ is isomorphic to $(J^{op},\Gamma_J^{op})$.
		
	\end{proof}
	
	\begin{df}
		If we have a group $G$ and an $n$-uple $\gamma=(g_1,...,g_n)\in G^n$, we denote $\gamma^{-1}=(g_1^{-1},...,g_n^{-1})$. 
	\end{df}
	\begin{df}
		For a group $G$, a $\ZZ\times G$-graded division algebra $\gamma_0\in G^n$ and $\gamma_1\in G^{m}$, we denote \[\cM(G,D,\gamma_0,\gamma_1)^{op}=(\cM(D)^{op},\Gamma(G,D,\gamma_0,\gamma_1)^{op}).\]
	\end{df}
	
	\begin{lemma}
		Let $D$ be a $(\ZZ\times G)$-graded division algebra, $\gamma_0=(g_1,...,g_{k_0})\in G^{k_0}$ and $\gamma_1=(g_{k_0+1},...,g_{k_0+k_1})\in G^{k_1}$. Then, $\cM(G,D,\gamma_0,\gamma_1)^{op}$ is isomorphic to the graded algebra $\cM(G,D^{op},\gamma_0^{-1},\gamma_1^{-1})$.
	\end{lemma}
	\begin{proof}
		Denote \[\gamma=((\delta_1,g_1),...,(\delta_{k_0+k_1},g_{k_0+k_1}))\] where $\delta_i=0$ if $0<i\leq k_0$ and $\delta_i=1$ otherwise. Then, if we denote $m=k_0+k_1$, $\cM(G,D,\gamma_0,\gamma_1)$ is the $\ZZ\times G$-graded algebra $\cM_m(D)$ with degrees given by $\deg (dE_{i,j})=(\delta_i-\delta_j,g_itg_j^{-1})$ for every $d\in D$ of degree $t$ and $0< i,j\leq m$. Denote:
		\[\gamma^{op}=((-\delta_1,g_1),...,(-\delta_m,g_m)).\]
		Now, $(\cM(D)^{op},\Gamma(G,D,\gamma^{op}))$ is the  $\ZZ\times G$-graded algebra $\cM_m(D^{op})$ with degrees given by $\deg (dE_{i,j})=(\delta_j-\delta_j,g_itg_j^{-1})$, which is the same as the algebra $\cM(G,D,\gamma_0,\gamma_1)^{op}$. 
		
		Consider the algebra homomorphism:

		\begin{equation*}
			\begin{split}
				\psi\colon \cM_m(D)^{op}&\to \cM_m(D^{op})\\
				d E_{i,j}&\mapsto dE_{j,i}.
			\end{split}
		\end{equation*}
		
		This is clearly an isomorphism of algebras and the grading induced in $\cM_m(D^{op})$ (i.e, the grading given by $\deg(\psi(x))=\deg(x)$) is $\Gamma(G,D^{op},(\gamma^{op})^{-1})$. Thus, \[(\cM(D)^{op},\Gamma(G,D,\gamma^{op}))\cong(\cM(D^{op}),\Gamma(G,D^{op},(\gamma^{op})^{-1})).\]
		Finally, since we have that $(\gamma^{op})^{-1}=((\delta_1,g_1^{-1}),...,(\delta_m g_m^{-1}))$,
		\[(\cM(D^{op}),\Gamma(G,D^{op},(\gamma^{op})^{-1}))\cong(\cM(G,D^{op},\gamma_0^{-1},\gamma_1^{-1}).\]
	\end{proof}
	
	We can summarize the results as follows:
	
	\begin{theorem}\label{t:NonGradedSimple}
		Let $(A,\varphi,\Gamma)$ be an $\Omega$-algebra satisfying $(T1)-(T4)$ and such that $(A,\Gamma)$ is not graded simple. Then, there is a $(\ZZ\times G)$-graded division simple algebra $D$, $k_0,k_1>0$, $\gamma_0\in G^{k_0}$ and $\gamma_1\in G^{k_1}$ consisting on distinct elements modulo $T=\supp (D)$, $\kappa_0\in \ZZ_{>0}^{k_0}$ and $\kappa_1\in \ZZ_{>0}^{k_1}$ such that $(A,\varphi,\Gamma)$ is isomorphic to $\cM(G,D,\kappa_0,\kappa_1,\gamma_0,\gamma_1)^{\ex}$.
		
		Moreover, $\cM(G,D,\kappa_0,\kappa_1,\gamma_0,\gamma_1)^{\ex}$ is isomorphic to $\cM(G,D,\kappa'_0,\kappa'_1,\gamma'_0,\gamma'_1)^{\ex}$ if either:
		\begin{itemize}
			\item[(1)] $D$ is isomorphic to $D'$ as graded algebras and there is $g\in G$ such that $\Xi(\kappa_i,\gamma_i)=g\Xi(\kappa'_i,\gamma'_i)$ for $i=0,1$.
			\item[(2)] $D$ is isomorphic to $D'^{op}$ as graded algebras and there is $g\in G$ such that $\Xi(\kappa_i,\gamma_i)=g\Xi(\kappa'_i,\gamma'^{-1}_i)$ for $i=0,1$.
		\end{itemize}
	\end{theorem}
	
	\subsection{Case with $(A,\Gamma)$ simple}
	
	Let $(A,\varphi,\Gamma)$ be an $\Omega$-algebra satisfying $(T1)-(T4)$ and such that $(A,\Gamma)$ is graded simple. Due to \cite[Theorem 2.6]{EKmon}, there is a $(\ZZ\times G)$-graded division algebra $D$ and a $(\ZZ\times G)$-graded right $D$-module $V$, such that $(A,\Gamma)$ is isomorphic to $\End_{D}(V)$ as a graded algebra. We will denote by $\Delta$ the grading of $V$. If we have a group homomorphism $\alpha\colon G\to H$, we are going to denote:
	
	\[^{\alpha}V=(V,\,^{\alpha}\Delta).\] 
	As in the previous section, we can assume that $\supp(^{\pi_1}V)=\{0,1\}$.
	
	\begin{remark}
		Let $G$ be an abelian group, $D$ and $D'$ two $G$-graded division algebras, a $G$-graded right $D$-module $V$ and a $G$-graded right $D'$-module $V'$. Given $g\in G$, an isomorphism of graded algebras $\psi_0\colon D\to D'$ and an isomorphism of graded modules $\psi_1\colon V\to V$ such that $\psi_1(vd)=\psi_1(v)\psi_0(d)$ for all $v\in V$ and $d\in D$, there is a unique homomorphism of graded algebras $\psi\colon \End_{D}(V)\to \End_{D'}(V')$ such that $\psi_1(rv)=\psi(r)\psi_1(v)$ for all $r\in \End_{D}(V)$ and $v\in V$. Moreover, all isomorphisms between these two graded algebras can be determined by a choice of a pair $(\psi_0,\psi_1)$ (\cite[Theorem 2.10]{EKmon}). In this situation, we denote $\psi=(\psi_0,\psi_1)$.
	\end{remark} 
	
	We are going to denote:
	\begin{equation}\label{eq:tau}\tau\colon \ZZ\times G\to \ZZ_2\times G
	\end{equation}
	the projection $\tau(i,g)=(\bar{i},g)$, where $\bar{i}$ denotes $i$ modulo $2$. Now, it follows that
	\[^{\tau}\End_{D}(V)=\End_{\,^{\tau}D}(^{\tau}V).\]
	If we denote $T=\supp(D)$, since $T$ is a finite group, as in the previous section, it is contained in $\{0\}\times G$. Thus, we can identifying with a subgroup of $G$ in a natural way, and we can write $D$ instead of $^{\tau}D$.
	
	\begin{remark}Any involution $\varphi$ of $\End_{D}(V)$ satisfying that $\deg(\varphi(x))=(-i,x)$ wherever $\deg(x)=(i,x)$ is an involution of the graded algebra $^{\tau}\End_{D}(V)$. 
	\end{remark}
	
	\begin{df} Given an abelian group $G$, a $G$-graded division algebra  $D$, an involution $\varphi_0$ of the graded algebra $D$, a $G$-graded right $D$-module $V$ and a nondegenerate $\varphi_0$-sesquilinear form $B\colon V\times V\to D$, we denote by $(\varphi_0,B)$ the involuion $\varphi$ of $\End_{D}(V)$ which satisfies:
		\[B(rv,w)=B(v,\varphi(r)w)\]
		for all $v,w\in V$ and $r\in \End_{D}(V)$. 
		
		If $\Gamma$ is the grading on $\End_{D}(V)$ induced by the grading on $D$ and the grading on $V$, we denote  as $\End(G,D,V,\varphi_0,B)$ the $G$-graded algebra with involution $(\End_{D}(V),(\varphi_0,B),\Gamma)$
	\end{df}
	
	\begin{proposition}\label{p:involutions}
		Let $(A,\varphi,\Gamma)$ be a triple satisfying $(T1)-(T4)$ such that $(A,\Gamma)$ is graded simple. Then, there exist a $(\ZZ\times G)$-graded-division algebra $D$, a graded right $D$-module $V$ with support in $\{0,1\}\times G$, an involution $\varphi_0$ of the graded algebra $D$, such that $(D,\varphi_0)$ is a simple algebra with involution and a nondegenerate hermitian or skew-hermitian $\varphi_0$-sesquilinear form $B\colon V\times V\to D$ such that the induced form $B\colon ^{\tau}V\times \,^{\tau}V\to D$ is graded of degree $(\overline{0},g)$ for some $g\in G$, such that $(A,\varphi,\Gamma)$ is isomorphic to $\End(\ZZ\times G,D,V,\varphi_0,B)$. Any such other pair $(\varphi_0',B')$ inducing the same involution on $\End_{D}(V)$ is of the form $(\mathrm{Int}(d)\circ \varphi_0,dB)$ where $d$ is a homogeneous element such that $\varphi_0(d)=\pm d$ and $\mathrm{Int}(d)(x)=dxd^{-1}$ for all $x\in D$.
		
		In addition, for every $5$-tuple $(G,D,V,\varphi_0,B)$ consisting on a group $G$, a $(\ZZ\times G)$-graded-division algebra $D$, a graded right $D$-module $V$ such that $\supp \,^{\pi_1}V=\{0,1\}\times G$, an involution $\varphi_0$ of the graded algebra $D$, such that $(D,\varphi_0)$ is a simple algebra with involution and a nondegenerate hermitian or skew-hermitian $\varphi_0$-sesquilinear form $B\colon V\times V\to D$ such that the induced form $B\colon ^{\tau}V\times \,^{\tau}V\to D$ is graded of degree $(\overline{0},g)$ for some $g\in G$, the graded algebra with involution $\End(\ZZ\times G, D,V,\varphi_0,B)$ satisfies $(T1)-(T4)$.
	\end{proposition}
	\begin{proof}
		As we showed earlier, $(A,\Gamma)$ is isomorphic to $\End_{D}(V)$ for a $\ZZ\times G$-graded division algebra $D$ and a $\ZZ\times G$-graded right $D$-module $V$. Moreover, we can assume that $\supp(^{\pi_1}V)=\{0,1\}$. Since $\varphi$ is an involution of $\End_{D}(^{\tau}V)$, \cite[Theorem 3.1]{EKR22} implies that there is an involution $\varphi_0$ of the graded algebra $D$ and a graded nondegenerate hermitian or skew-hermitian $\varphi_0$-sesquilinear form $B\colon ^{\tau}V\times^{\tau}V\to D$ satisfying that $B(rv,w)=B(v,\varphi(r)w)$ for all $v,w\in V$ and $r\in \End_{D}(V)$.
		
		Assume that the degree of $B$ is $(\bar{i},g)\in\ZZ_2\times G$. In order to show that $\overline{i}=\subo$, we are going to assume that $\bar{i}=\subuno$ and argue by contradiction.
		
		Take an element $h_1\in G$ such that $^{\tau}V_{\subo,h_1}\neq 0$ and a nonzero element $v\in \,^{\tau}V_{\subo,h_1}$. Clearly, since $\supp (^{\pi_1}V)= \{0,1\}$, it follows that $v\in V_{0,h_1}$. The fact that $B$ is nondegenerate, $\supp(D)\subseteq\{\subo\}\times G$ and it has degree $(\subuno,g)$, implies that there is $h_2\in G$ and $w\in \,^{\tau}V_{(\subuno,h_2)}$ such that $B(v,w)\neq 0$. As argued with $v$, we have that  $w\in V_{(1,h_2)}$. Let $\{v,w,u_1,...,u_k\}$ be a homogeneous basis of $V$ as a $D$ module and let $r\in \End_{D}(V)$ be the homomorphism such that $r(w)=v$, $r(v)=0$ and $r(u_i)=0$ for all $i\in \{1,...,k\}$. $r$ is a homomorphism of degree $(-1,h_2h_1^{-1})$. Therefore, $\varphi(r)$ is a homomorphism of degree $(1,h_2h_1^{-1})$. Hence, $\varphi(r)w=0$. This leads to a contradiction since:
		
		\[0\neq B(v,w)=B(r(w),w)=B(w,\varphi(r)w)=0.\]
		
		If we have another pair $(\varphi'_0, B')$ inducing the same involution on $\End_{D}(V)$, \cite[Theorem 3.1]{EKR22} implies that there is a homogeneous element $d\in D$ such that $\varphi_0(d)=\pm d$ and such that $(\varphi'_0,B')=(\mathrm{Int}(d)\circ \varphi_0,dB)$.`
		
		To prove the second part of the proposition, we take $D,\varphi_0, V$ and $B$ as in the statement. The fact that $\End_{D}(V)$ is a simple algebra with involution implies that $(\End_{D}(V),\varphi)$ is a simple associative algebra with involution, which implies $T1$. Since $D$ and $V$ are $(\ZZ\times G)$-graded, then (T2) follows. Since $\supp(^{\pi_1}D)=\{0\}$ and $\supp(^{\pi_1}V)=\{0,1\}$ (T3) follows. 
		
		Finally, we need to prove (T4). We denote by $\varphi$ the involution $(\varphi_0,B)$. If $r\in \End_{D}(V)_{(0,g)}$ for some $g\in G$. Since $\varphi$ is an involution of $^{\tau}\End_{D}(V)$, it follows that $\varphi(r)\in \End_{D}(V)_{(0,g)}$. If $r\in \End_{D}(V)_{(1,h)}$ for some $h\in G$. Since $\varphi$ is an involution of $^{\tau}\End_{D}(V)$, it follows that there are $r_i\in \End_{D}(V)_{(i,h)}$ with $i=\pm 1$ such that $\varphi(r)=r_{-1}+r_{1}$. We just need to show that $r_1=0$. We will argue by contradiction.
		
		Assume that $r_1\neq 0$. Then, there is $h_1\in G$ and $v\in V_{(0,h_1)}$ such that $0\neq w=r_1(v)\in V$. Moreover, $r_{-1}(v)=0$ due to the fact that $\supp( ^{\pi_1}V)=\{0,1\}$. Therefore, $\varphi(v)=w$. Since $B$ is a graded form of degree $(\subo,g)$ of $^{\tau}V$, there is an element $h_2\in G$ and $u\in V_{(1,h_2)}$ such that $B(u,v)\neq 0$. Since $r\in \,^{\tau} \End_D(V)_{(1,h)}$, it follows that $r(u)=0$.  Thus:
		
		\[0\neq B(u,w)=B(u,\varphi(r)(v))=B(r(u),v)=0,\]
		which is a contradiction.
		
	\end{proof}
	
	\section{Gradings on the associative algebras over algebraically closed fields}\label{sec:AlgClosedField}
	
	In this section we will assume that $\FF$ is an algebraically closed field of characteristic other than $2$ and $G$ an abelian group.
	
	\subsection{Graded-division algebras.}
	
	The first step to specializing the results of the previous chapter to algebraically closed fields is to classify $G$-graded-division simple algebras with involution. We will start recalling the case when the underlying algebra is simple. The arguments can be found in \cite[Chapter 2]{EKmon}
	
	Let $T$ be a finite subgroup of $G$ such that $\mathrm{char}(\FF)$ doesn't divide its order and we let $\beta\colon T\times T\to\FF$ be a nondegenerate alternating bicharacter on T ,i.e, a map which is multiplicative on each variable and satisfies $\beta(t,t)=1$ for all $t\in T$. Due to the fact that $\beta$ is alternating and nondegenerate, we can decompose $T$ as:
	
	\[T=H_1'\times H_1''\times ...\times H_r'\times H''_r\]
	for some $r\in \ZZ_{>0}$ and some groups $H'_1,H''_1,...,H'_r,H''_r$ in such way that $H'_i\times H''_i$ and $H'_j\times H''_j$ are orthogonal for $i\neq j$ and such that $H'_i$ and $H''_i$ are cyclic groups of order $l_i$. Choose generators $a_i$ and $b_i$ of $H'_i$ and $H''_i$ respectively and denote $\epsilon_i=\beta(a_i,b_i)$. We can endow $\cM_{l_1}(\FF)\otimes_{\FF}...\otimes_{\FF}\cM_{l_r}(\FF)$ with a $G$-grading given by:
	
	\[\deg(I\otimes...\otimes I\otimes X_i\otimes I\otimes ...\otimes I)=a_i\]
	and
	\[\deg(I\otimes...\otimes I\otimes Y_i\otimes I\otimes ...\otimes I)=b_i\]
	where:
	\[X_i=\sum_{k=1}^{l_i}\epsilon_i^{k-1}E_{k,k},\text{ and } Y_i=\sum_{k=1}^{l_i}E_{k,k+1}\]
	with the indices modulo $l_i$ for each $i\in \{1,...,r\}$. With this grading, the algebra $\cM_{l_1}(\FF)\otimes_{\FF}...\otimes_{\FF}\cM_{l_r}(\FF)$ is a graded division algebra.  
	
	\begin{df}
		The previous grading, induces a $G$-grading naturally on $\cM_{l_1\cdots l_r}(\FF)$ via the Kronecker product. We call this graded division algebra a \emph{standard realization} of $T$ and $\beta$ which we denote $\cD(T,\beta)$.
	\end{df}
	
	\begin{remark}\label{rm:Xt1Xt2Commuting}
		With te previous notation, we denote $X_{a_i}=I\otimes...\otimes I\otimes X_i\otimes I\otimes ...\otimes I$, $X_{b_i}=I\otimes...\otimes I\otimes Y_i\otimes I\otimes ...\otimes I$ and $X_{(a_1^{i_1},b_1^{j_1},...,a_r^{i_r},b_r^{j_r})}=X_{a_1}^{i_1}X_{b_1}^{j_1}...X_{a_r}^{i_r}X_{b_r}^{j_r}$.
		
		We can realise that $X_{a_i}^{l_i}=1=X_{b_i}^{l_1}$  for each $i\in\{1,...,r\}$ and also that $X_{t_1}X_{t_2}=\beta(t_1,t_2)X_{t_2}X_{t_1}$ for every $t_1,t_2\in T$.
	\end{remark}
	The standard realizations of $T$ and $\beta$ are all isomorphic. Therefore, the notation $D(T,\beta)$, is a well defined notation up to isomorphism. The following Theorem classifies all the finite dimensional $G$-graded division algebra
	\begin{theorem}[\cite{EKmon}]\label{t:DivAlgClass} Let $T$ be an finite abelian group and let $\FF$ be an algebraically closed field. There exist a grading on the matrix algebra $\cM_n(\FF)$ with support $T$ marking $\cM_n(\FF)$ a graded-division algebra if and only if $\mathrm{char}(\FF)$ does not divide $n$ and $T\cong \ZZ_{l_1}^2\times \dots\ZZ_{l_r}^2$ with $l_1\cdots l_r=n$. The isomorphism classes of such gradings are in one to one correspondence with the isomorphism classes of nondegenerate alernating bicharacters $\beta\colon T\times T\to \FF$.
		
	\end{theorem}

	In order to classify algebras up to isomorphism, we will need the following result:
	
	\begin{lemma}\label{l:DivAlgOp}
		Let $T\cong \ZZ^2_{l_1}\times\dots\ZZ^2_{l_r}$ be a subgroup of $\ZZ\times G$ and $\mathrm{char}(\FF)$ not dividing $l_1\cdots l_r=n$, and a nondegenerate alternating bicharacter $\beta\colon T\times T\to G$. Then, the $\ZZ\times G-D(T,\beta)$ satisfies that $D(T,\beta)^{op}\cong D(T,\beta\circ \ex)$ where $\ex\colon T\times T\to T\times T$ is the map given by $\ex(t_1,t_2)=(t_2,t_1)$.
	\end{lemma}
	
	\begin{proof}
		This is due to the fact that $T\subseteq \{0\}\times G$ and due to the fact that as shown in Remark \ref{rm:Xt1Xt2Commuting}, $X_{t_1}X_{t_2}=\beta(t_1,t_2)X_{t_2}X_{t_1}$.
	\end{proof}
	
	\begin{remark}\label{rm:bicharacterOp}
		For a  nondegenerate alternating bicharacter $\beta\colon \ZZ_2^n\times \ZZ_2^n\to \FF$, since $\beta$ is multiplicative, its image is is $\{\pm 1\}$. Thus, since $\beta$ is alternating, \[\beta(u,v)=\beta(v,u)^{-1}=\beta(v,u)\] which implies that $\beta=\beta\circ \ex$.
	\end{remark}
	
	With that in mind, we can star the classification of the $G$-graded-simple graded-division simple algebras with involution. To begin with, we will approach the case when the algebra is simple.
	
	\begin{lemma}\label{l:2group1Dim}
		Let $(D,\varphi_0)$ be a $G$-graded-simple graded-division simple algebra with involution over an algebraically closed field $\FF$. Then, $T=\supp D$ is an elementary abelian $2$-group. Moreover, every homogeneous component is one-dimensional.
	\end{lemma}
	\begin{proof}
		Since $D_e$ is a division algebra and $\FF$ is an algebraically closed field, it implies that $D_e$ is one-dimensional. Thus, due to \cite[Lemma 3.5]{EKR22}, $T$ is an elementary abelian $2$-group. Now, if $t,t'\in T$ and $0\neq X_t X_{t'}\in D_t$, it follows that  $X_tX_{t'}^{-1}\in D_e$. Thus, $D_t=X_tD_e$, which implies that every homogeneous component is one-dimensional.
	\end{proof}
	
	\begin{df}\label{df:Dinv}
		Let $T$ be an elementary abelian $2$-subgroup of $\ZZ\times G$, $\beta\colon T\times T\to G$ a nondegenerate alternating bicharacter and $\tau\colon T\to \FF$ a quadratic form on $T$ whose polar form is $\beta$, i.e., such that $\beta(t_1,t_2)=\tau(t_1t_2)\tau(t_1)\tau(t_2)$. We denote by $D_{\mathrm{inv}}(T,\beta,\tau)$ the graded algebra $D(T,\beta)$ with the involution given by:
		\[d\mapsto \tau(t)d.\]
		On the graded algebra $D(T,\beta)$, the transposition, is an involution. We also denote by $\beta$, the quadratic form for which the involution is given by:
		\[d\mapsto \beta(t)d\]
		for every $t\in T$ and $d\in D_t$.
	\end{df}
	
	\begin{proposition}\label{p:GDivAlgWithInv}
		Let $T$ be an abelian group, $\FF$ and algebraically closed field of characteristic different from $2$ and $\beta\colon T\times T\to \FF$ a nondegenerate alternating bicharacter. Then, $D(T,\beta)$ admits an involution if and only if $T$ is an elementary $2$-group. Then, this algebra with involution is of the form $\Dinv(T,\beta,\tau)$ for a quadratic form $\tau$ whose polar form is $\beta$. Any involution of $D(T,\beta)$ is of the form:
		
		\[X\mapsto X_t^{-1}X^{t}X_t\]
		for a uniquely determined $t\in T$ and an element $X_t\in D_t$.
	\end{proposition}
	\begin{proof}
		This proof is in \cite[Proposition 2.51]{EKmon}
	\end{proof}
	
	\begin{remark}\label{rm:StandardRealizationTrans}
		We denote $\Dinv(T,\beta)$ the graded algebra $D(T,\beta)$ with the involution given by the transposition.
	\end{remark}
	\begin{df}
		Let $(A,\varphi,\Gamma)$ be a $G$-graded algebra with involution and let $t\in G\setminus\supp(\Gamma)$ be an order $2$ element. We denote $(A,\varphi,\Gamma)^{\ex}_t$ the $G$-graded algebra $(A\oplus A^{op},\ex,\Gamma^{ex}_t)$ where $\Gamma^{\ex}_t$ is a $G$-grading given by $\deg(x,\varphi(x))=h$ and $\deg(x,-\varphi(x))=ht$ if $x\in A_h$.
	\end{df}
	
	The fact that the previous grading is well defined is shown in \cite{BGGradedInvol}.
	
	\begin{df}
		Let $T$ be an elementary abelian $2$-group and $T'$ be a subgroup of index $2$, let $t\in T\setminus T'$ and let $\beta\colon T'\times T'\to \FF$ be an alternating bicharacter. We denote $\beta^{[t]}\colon T\times T\to \FF$, the bicharacter given by:\[\beta^{[t]}(ut^{k_1},vt^{k_2})=\beta(u,v)\]
		for every $u,v\in T'$ and $k_1,k_2\in \{0,1\}$. Similarly, if $\tau\colon T'\to \FF$ is a quadratic form, we define the quadratic form $\tau^{[t]}\colon T\to \FF$ as the map given by:
		\[\tau^{[t]}(ut^{k})=\tau(u)(-1)^k\]
		for every $u\in T'$ and $k\in\{0,1\}$.
	\end{df}
	
	We can notice that if $\beta$ is the polar form of $\tau$, then, $\beta^{[t]}$ is the polar form of $\tau^{[t]}$.
	
	\begin{remark}\label{rm:involTau}Let $T$ be an elementary abelian $2$-group and $T'$ be a subgroup of index $2$, let $t\in T\setminus T'$, let $\beta\colon T'\times T'\to \FF$ be an alternating bicharacter and let $\tau\colon T'\to \FF$ be a quadratic form whose polar form is $\beta$. The involution on $\Dinv(T',\beta,\tau)^{\ex}_t$ is given by:
		\[d\mapsto \tau^{[t]}(s)d\]
		for every homogeneous $d$ of degree $s$. For the case of $\Dinv(T',\beta)^{\ex}_t$, this implies that the involution is given by:
		\[(X,Y)\mapsto (Y^t,X^t).\]
	\end{remark}
	
	\begin{proposition}\label{p:ExGrDivAlg}
		Let $T$ be an elementary abelian $2$-subgroup of $G$ of rank $2m+1$ for some $m\in \ZZ_{\geq0}$, let $t\in T$ be an order $2$ element and let $T_1$ and $T_2$ be two elementary $2$-subgroups of index $2$ not containing $t$. Let $\tau_1\colon T_1\to \FF$ and $\tau_2\colon T_2\to \FF$ be two quadratic forms such that their polar forms, $\beta_1\colon T_1\times T_1\to\FF$ and $\beta_2\colon T_2\times T_2\to\FF$, are two nondegenerate alternating bicharacters and such that $\tau_1^{[t]}=\tau_2^{[t]}$. Then $\Dinv(T_1,\beta_1,\tau_1)^{\ex}_t\cong \Dinv(T_2,\beta_2,\tau_2)^{\ex}_t$.
	\end{proposition}
	\begin{proof}
		We will denote by $(D,\varphi_0)$ the underlying algebra with involution from the graded algebra with involution  $\Dinv(T_1,\beta_1,\tau_1)$, and we write:
		\[\Dinv(T_1,\beta_1,\tau_1)_{T_2}=\bigoplus_{h\in T_2}(\Dinv(T_1,\beta_1,\tau_1)^{\ex}_t)_h.\]
		Then, $\Dinv(T_1,\beta_1,\tau_1)_{T_2}$ is a $G$-graded-division algebra and the exchange involution restricts to $\Dinv(T_1,\beta_1,\tau_1)_{T_2}$. Let \[\psi\colon \Dinv(T_1,\beta_1,\tau_1)_{T_2}\to D\]
		be the projection $\psi(x,y)=x$. We claim that $\psi$ is an algebra isomorphism. Indeed, it is surjective because if $x$ is a homogeneous element of degree $h$ in $D$, either $h\in T_2$ or $ht\in T_2$. Thus, either $(x,\varphi(x))\in \Dinv(T_1,\beta_1,\tau_1)_{T_2}$ or $(x,-\varphi(x))\Dinv(T_1,\beta_1,\tau_1)_{T_2}$. In order to show that it is an isomorphism, we use that Lemma \ref{l:2group1Dim} implies that:
		\[\dim(\Dinv(T_1,\beta_1,\tau_1)_{T_2})=\lvert T_2\rvert  =\lvert T_1\rvert=\dim(D)\]
		where $\lvert H\rvert$ means the cardinal of $H$. $\psi$ induces a grading and an involution on $D$ which we denote $\Gamma_{T_2}$ and $\bar{\varphi}_0$.
		
		Since $T_2$ and $T_1$ have index $2$ on $T$, there is a unique basis ${Y_h}_{h\in T_2}$ of $D$ defined by:
		\[Y_h=\psi(X_{ht^k},(-1)^{k}\varphi_0(X_{ht^k}))\]
		for a $k\in\{0,1\}$ such that $ht^k\in T_1$. Clearly, this basis is a basis as the one in Remark \ref{rm:Xt1Xt2Commuting}.
		
		We are going to prove that $(D,\bar{\varphi}_0,\Gamma_{T_2})$ is $D(T_2,\beta_2,\tau_2)$. In order to do so, we need to show that $Y_{h_1}Y_{h_2}=\beta_2(h_1,h_2)Y_{h_2}Y_{h_1}$ and we need to show that $\bar{\varphi}_0(Y_h)=\tau_2(h)Y_h$ for all $h_1,h_2,h\in T_2$.
		
		Given $h_1,h_2\in T_2$, let $k_1,k_2\in \{0,1\}$ be such that $h_1t^{k_1},h_2t^{k_2}\in T_1$. Then:
		\begin{equation*}
			\begin{split}
				Y_{h_1}Y_{h_{2}}&=\psi(X_{h_1t^{k_1}},(-1)^k\varphi_0(X_{h_1t^{k_1}}))\psi(X_{h_2t^{k_2}},(-1)^k\varphi_0(X_{h_2t^{k_2}}))\\
				&=\beta_1(h_1t^{k_1},h_2t^{k_2})\psi(X_{h_2t^{k_2}},(-1)^k_2\varphi_0(X_{h_2t^{k_2}}))\psi(X_{h_1t^{k_1}},(-1)^k\varphi_0(X_{h_1t^{k_1}}))\\
				&=\beta_1(h_1t^{k_1},h_2t^{k_2})Y_{h_{2}}Y_{h_1}
			\end{split}
		\end{equation*}
		and after using that 
		\[\beta_1(h_1t^{k_1},h_2t^{k_2})=\beta_1^{[t]}(h_1t^{k_1},h_2t^{k_2})=\beta_2^{[t]}(h_1t^{k_1},h_2t^{k_2})=\beta_2(h_1,h_2),\]
		we get 
		\[Y_{h_1}Y_{h_2}=\beta_2(h_1,h_2)Y_{h_2}Y_{h_1}.\]
		
		Similarly, for $h\in T_2$ and $k\in\{0,1\}$ such that $ht^k\in T_1$, using Remark \ref{rm:involTau} and the fact that $(X_{ht^k},(-1)^k\varphi_0(X_{ht^k}))$ has degree $h$, it follows that 
		
		\begin{equation*}
			\begin{split}
				\bar{\varphi_0}(Y_h)&=\bar{\varphi_0}(\psi(X_{ht^{k}},(-1)^k\varphi_0(X_{ht^{k}}))\\
				&=\psi(\ex(X_{ht^{k}},(-1)^k\varphi_0(X_{ht^{k}}))\\
				&=\tau_1(h)\psi(X_{ht^{k}},(-1)^k\varphi_0(X_{ht^{k}})\\
				&=\tau_1(h)Y_h,
			\end{split}
		\end{equation*}
		which due to the fact that 
		\[\tau_1(h)=\tau^{[t]}_1(ht^k)(-1)^k=\tau_2^{[t]}(ht^k)(-1)^k=\tau_2(h)(-1)^k(-1)^k=\tau_2(h),\]
		implies:
		\[\bar{\varphi}_0(Y_h)=\tau_2(h)Y_h.\]
		Therefore, $(D,\bar{\varphi}_0,\Gamma_{T_2})=\Dinv(T_2,\beta_2,\tau_2)$.
		
		Finally, due to the fact that $\psi$ is a homogenous isomorphism and the fact that $(1,-1)$ has degree $t$ in $\Dinv(T_1,\beta_1,\tau_1)^{\ex}_t$, implies that if $x$ is homogeneous of degree $h$ in $\Dinv(T_2,\beta_2,\tau_2)$, then $(x,\varphi(x))$ has degree $h$ and $(x,-\varphi(x))$ has degree $ht$ in $\Dinv(T_2,\beta_2,\tau_2)^{\ex}_t$, implying that $\Dinv(T_1,\beta_1,\tau_1)^{\ex}_t=\Dinv(T_2,\beta_2,\tau_2)^{\ex}_t$.
		
	\end{proof}
	
	\begin{theorem}\label{t:NonSimpleGradedDivAlg}
		Let $T$ be an elementary $2$-subgroup of rank $2m$ of $G$ and let $t\in G\setminus T$ be an order $2$ element, $\tau\colon T\to \FF$ a quadratic form such that its polar form $\beta\colon T\times T\to \FF$  is a nondegenerate alternating bicharacter. Then $\Dinv(T,\beta,\tau)^{\ex}_t$ is a graded-division simple algebra with involution such that the underlying algebra is not simple. Moreover, given  a $G$-graded-division simple algebra with involution $(D,\varphi_0,\Gamma)$ with $D$ not simple, there is an elementary $2$-subgroup $T$ of rank $2m$ of $G$, an order $2$ element $t\in G\setminus T$ and  a quadratic form $\tau\colon T\to \FF$ whose polar form is a nondegenerate alternating bicharacter $\beta\colon T\times T\to \FF$ such that \[(D,\varphi_0,\Gamma)\cong \Dinv(T,\beta,\tau)^{\ex}_t.\] 
		
		The algebras $\Dinv(T_1,\beta_1,\tau_1)^{\ex}_{t_1}$ and $\Dinv(T_2,\beta_2,\tau_2)^{\ex}_{t_2}$ are isomorphic if and only if $T_1\langle t_1\rangle=T_2\langle t_2\rangle$, $t_1=t_2$, $\beta_1^{[t_1]}=\beta_2^{[t_2]}$ and $\tau_1^{[t_1]}=\tau_2^{[t_2]}$.
	\end{theorem}
	\begin{proof}
		Denote $(D,\varphi_0,\Gamma)= \Dinv(T,\beta,\tau)$. The homogeneous component of the graded algebra with involution $\Dinv(T,\beta,\tau)$ have dimension $1$. Thus, every homogeneous component of $\Dinv(T,\beta,\tau)^{\ex}_t$ has dimension 1, which implies that $(x,y)$ is homogeneous in $\Dinv(T,\beta,\tau)^{\ex}_t$, if and only if $x$ is homogeneous in $\Dinv(T,\beta,\tau)$ and $y=\pm\varphi_0(x)$. Since $\Dinv(T,\beta,\tau)$ is a graded division algebra, it follows that $x$ and $y$ are invertible and as a consequence, $(x,y)$ is invertible. Therefore, $\Dinv(T,\beta,\tau)^{\ex}_t$ is a graded division algebra. Moreover, the fact that $D$ is simple, implies that $(D,\varphi_0)$ is a simple algebra with involution.
		
		For the second part, take a $G$-graded-division simple algebra with involution $(D,\varphi_0,\Gamma)$ for which $D$ is not simple. Since $D$ is not simple, there is a proper ideal $0\neq I\subseteq D$ which satisfies:
		\[(D,\varphi_0)\cong (I\oplus I^{op},\ex).\]
		For the proof we are going to identify both algebras with involution. The fact that $(D,\Gamma)$ is an algebra with involution implies that its support is a group, so we can assume that the grading group $T$ is the support. The center of $D$:
		\[\cZ(D)=\FF(1,1)\oplus \FF(1,-1)\]
		is a graded subalgebra and the involution restricts to it. Therefore, the space of symmetric elements and the space of antisymmetric elements are graded, so $\deg(1,1)=e$ and $\deg(1,-1)=t$ for an order $2$ element $t\in T$. Denote by $\overline{\Gamma}$, the coarsening by the group $\overline{T}=T/\langle t\rangle$ given by:
		\[\overline{\Gamma}\colon D=\bigoplus_{\bar{t}'\in\overline{T} }D_{\bar{t}'}\]
		where $D_{\bar{t}'}=D_{t'}\oplus D_{t't}$ (recall the definition of coarsening from Section \ref{sec:Prelim}). Lemma \ref{l:2group1Dim} implies that $T$ is a $2$-group and due to this fact, we can identify $\overline{T}$ with a subgroup $T'$ of $T$ in such way that $T\cong T'\times \langle t\rangle$. In this coarsening, $\cZ(D)_e=\cZ(D)$, which implies that $e_1=(1,0)$ and $e_2=\varphi_{0}(e_1)=(0,1)$ are homogeneous elements of degree $e$. This implies that the subspaces $I=De_1$ and $\varphi_0(I)=De_2$ are graded, and for every $x\in I_{\bar{t'}}$, $\varphi_0(x)\in \varphi_0(I)_{\bar{t'}}$. We denote by $\Gamma_I$ the restriction of the grading to $I$. Since the homogeneous components of $D$ in $\overline{\Gamma}$ have dimension $2$, it follows that $D_{\bar{t'}}=\FF x\oplus \FF\varphi_0(x)$, implying that the homogeneous components of $I$ have dimension $1$. Therefore, for every $x\in I_{\bar{t'}}$, there is a unique $y\in I$, concretely $y\in T_{\bar{t'}}$, such that $x+\varphi_0(y)\in D_{t'}$, and multiplying by $(1,-1)$, $x-\varphi_0(y)\in D_{t't}$. We can define an involution $\bar{\varphi}_0$ in $I$ by sending $x$ to $y$ and using the identification of $D$ and $I\oplus I^{op}$, the grading is given by \[D_{t'}=\{(x,\bar{\varphi_0}(x))\mid x\in I_{t'}\}\] and 
		\[D_{t't}=\{(x,-\bar{\varphi_0}(x))\mid x\in I_{t'}\}\]
		for every $t'\in T'$. We can show that $(I,\Gamma_I)$ is a graded algebra since for a homogeneous element $x\in I$, $(x,\bar{\varphi}_0)$ is a homogeneous element of $(D,\Gamma)$. Since $(D,\Gamma)$ is a graded division algebra, $(x,\bar{\varphi}_0)$ is an invertible element and this implies that $x$ is invertible. Since $(I,\bar{\varphi_0},\Gamma_I)$ is a graded division algebra with involution and $I$ is simple, Proposition \ref{p:GDivAlgWithInv} implies that there is a quadratic form $\tau\colon T\to \FF$ such that its polar form is the nondegenerate alternating bicharacter $\beta\colon T\times T\to \FF$, such that \[(I,\bar{\varphi}_0,\Gamma_I)\cong \Dinv(T,\beta,\tau).\]
		Hence:
		\[(D,\varphi_0,\Gamma)\cong \Dinv(T,\beta,\tau)^{\ex}_t.\] 
		
		We can now tackle the isomorphism problem. Assume first that $T_1\langle t_1\rangle=T_2\langle t_2\rangle$, $t_1=t_2$, $\beta_1^{[t_1]}=\beta_2^{[t_2]}$ and $\tau_1^{[t_1]}=\tau_2^{[t_2]}$ due to Proposition \ref{p:ExGrDivAlg} the result follows. Conversely, assume that the algebras  $\Dinv(T_1,\beta_1,\tau_1)^{\ex}_{t_1}$ and $\Dinv(T_2,\beta_2,\tau_2)^{\ex}_{t_2}$ are isomorphic.
		
		The fact that:
		\[\{e,t_1\}=\supp(\cZ(\Dinv(T_1,\beta_1,\tau_1)^{\ex}_{t_1}))=\supp \Dinv(\cZ((T_2,\beta_2,\tau_2)^{\ex}_{t_2}))=\{e,t_2\}\]
		implies that $t_1=t_2$. Moreover,
		\[T_1\langle t_1\rangle=\supp(\Dinv(T_1,\beta_1,\tau_1)^{\ex}_{t_1})=\supp(\Dinv(T_2,\beta_2,\tau_2)^{\ex}_{t_2})=T_2\langle t_2\rangle.\]
		Finally, given $h_1,h_2\in T_1$, $X_{h_1},X_{h_2}\in (\Dinv(T_1,\beta_1,\tau_1)$ and $k_1,h_2\in \{0,1\}$, the elements $(X_{h_1},(-1)^{k_1}\varphi_{0}X_{h_1})$ and $(X_{h_2},(-1)^{k_2}\varphi_{0}X_{h_2})$ are of degrees $h_1t^{k_1}$ and $h_2t^{k_2}$ respectively. Now due to Remark \ref{rm:bicharacterOp} the following identity holds:
		\begin{multline*}(X_{h_1},(-1)^{k_1}\varphi_0(X_{h_1}))(X_{h_2},(-1)^{k_2}\varphi_0(X_{h_2}))\\=\beta_1(h_1,h_2)(X_{h_2},(-1)^{k_2}\varphi_0(X_{h_2}))(X_{h_1},(-1)^{k_1}\varphi_0(X_{h_1})).
		\end{multline*}
		Thus, the fact that the homogeneous components have dimension $1$, implies that for elements $x,y$ of degree $h_1,h_2\in \Dinv(T_1,\beta_1,\tau_1)^{\ex}_{t_1}$, we have the following identity:
		\[xy=\beta_1^{[t_1]}(h_1,h_2)yx.\]
		Doing the same for $ \Dinv(T_2,\beta_2,\tau_2)^{\ex}_{t_2}$ it follows that $\beta_1^{[t_1]}=\beta_2^{[t_2]}$. In the same way, using Remark \ref{rm:involTau}, it follows that $\tau_1^{[t_1]}=\tau_2^{[t_2]}$.
	\end{proof}
	
	\begin{remark}
		Let $T$ be a $2$-elementary abelian subgroup of a group $G$ and $\beta\colon T\times T\to \FF$ a nondegenerate abelian bicharacter and $t\in G\setminus T$ an order $2$ element. With the notation of Remark \ref{rm:Xt1Xt2Commuting} and Remark \ref{rm:StandardRealizationTrans}, for every $t'\in T$ and $k\in\{0,1\}$, we denote
		\[Y_{t't^{k}}=(X_{t'},(-1)^{k}X_{t'}^t).\]
		Taking the quadratic form $\tau\colon T\to \FF$ such that $X_{t'}^t=\tau(t')X_{t'}$ for every $t'\in T$, clearly, for every $t'\in T\langle t\rangle$, it follows that $\varphi(Y_{t'})=\tau^{[t]}T_{t'}$.
	\end{remark}
	
	\begin{df}
		Given a group $G$, a $2$-elementary abelian subgroup of $G$ of rank $2m$, an order $2$ element $t\in G\setminus T$, let $\tau\colon T\to \FF$ be a quadratic form whose polar form is a nondegenerate alternating bicharacter $\beta\colon T\times T\to \FF$. Denote $(D,\varphi,\Gamma)=\Dinv(T,\beta,\tau)^{\ex}_t$. For a homogeneous element $d$, we denote $\Dinv(T,\beta,\tau,d)^{ex}_t=(D,\mathrm{Int}(d)\circ \varphi,\Gamma)$, where $\mathrm{Int}(d)(x)=dxd^{-1}$ for all $x\in D$.
	\end{df}
	
	\begin{proposition}\label{p:RemoveTau}
		Let $T_1$ and $T_2$ be two elementary $2$-subgroups of $G$. Let $t_i\in  G\setminus T_i$  for $i\in\{1,2\}$ be two order $2$ elements, let $\tau_1\colon T\to \FF$ and $\tau_2\colon T\to \FF$ be two quadratic form whose polar forms are respectively the nondegenerate alternating bicharacters $\beta_1\colon T\times T\to \FF$ and $\beta_2\colon T\times T\to \FF$. Then, there exists a homogeneous element $d\in \Dinv(T_1,\beta_1,\tau_1)^{\ex}_{t_1}$ such that $\Dinv(T_1,\beta_1,\tau_1,d)^{\ex}_{t_1}\cong \Dinv(T_2,\beta_2,\tau_2)^{\ex}_{t_2}$ if and only if $T_1\langle t_1\rangle=T_2\langle t_2\rangle$, $t_1=t_2$ and $\beta_1^{[t_1]}=\beta_2^{[t_2]}$.
	\end{proposition}
	
	\begin{proof}
		Assume that there exists a homogeneous element $d\in \Dinv(T_1,\beta_1,\tau_1)^{\ex}_{t_1}$ such that $\Dinv(T_1,\beta_1,\tau_1,d)^{\ex}_{t_1}\cong \Dinv(T_2,\beta_2,\tau_2)^{\ex}_{t_2}$. We can prove that  $T_1\langle t_1\rangle=T_2\langle t_2\rangle$, $t_1=t_2$ and $\beta_1^{[t_1]}=\beta_2^{[t_2]}$ in the same way as in Theorem \ref{t:NonSimpleGradedDivAlg}.
		
		Assume now that $T_1\langle t_1\rangle=T_2\langle t_2\rangle$, $t_1=t_2$ and $\beta_1^{[t_1]}=\beta_2^{[t_2]}$. Denote $\tau'\colon T_1\to \FF$ the map given by $\tau'(t_1)=\tau_2^{[t_2]}(t_1)$ for all $t_1\in T_1$. Theorem \ref{t:NonSimpleGradedDivAlg} implies that $\Dinv(T_2,\beta_2,\tau_2)^{\ex}_{t_2}\cong \Dinv(T_1,\beta_1,\tau')^{\ex}_{t_1}$.
		
		Denote $(D,\varphi_1,\Gamma_1)=\Dinv(T_1,\beta_1,\tau_1)^{\ex}_{t_1}$ and $(D,\varphi_2,\Gamma_2)=\Dinv(T_1,\beta_1,\tau')^{\ex}_{t_1}$.
		
		Consider the vector space endomorphism:
		
		\[\Psi\colon D\to D\]
		
		induced by $(x,0)\mapsto (x,0)$ and $(0,x)\mapsto (0,\tau_1(t)\tau'(t)x)$ for every homogeneous element $x\in D(T_1,\beta_1)$ of degree $t$. This is clearly an isomorphism of vector spaces. The fact that it is the identity in the first component and it is a composition of two involutions in the second component implies that it is an isomorphism of algebras. It sends the grading $\Gamma_2$ to the grading $\Gamma_1$ since it sends $(x,\tau'(t)x)$ to $(x,\tau_1(t)x)$ and $(x,-\tau'(t)x)$ to $(x,-\tau_1(t)x)$ for every $x\in D(T_1,\beta_1)$ of degree $t$. From this, it follows that:
		
		\[(D,\varphi_2,\Gamma_2)\cong(D,\Psi\circ \varphi_2\circ \Psi^{-1},\Gamma_1).\]
		
		The involution $\Psi\circ \varphi_2\circ \Psi^{-1}$ is given by:
		
		\[\Psi\circ \varphi_2\circ \Psi^{-1}(x,\pm\tau_1(t)x)=\tau'(t)(x,\pm\tau_1(t)x).\]
		for every $x\in D(T_1,\beta_1)$ of degree $t$. Finally using Proposition \ref{p:GDivAlgWithInv}, it follows that there is a uniquely determined $t'\in T_1$ such that $\tau'(t)x=\tau_1(t)X_{t'}xX_{t'}$ for all $x\in D(T_1,\beta_1)$ of degree $t$. And so, $\Psi\circ \varphi_2\circ \Psi^{-1}=\mathrm{Int}(Y_{t'})\circ \varphi_2$, and the statement follows.
		
	\end{proof}
	
	\subsection{Associative algebras}
	Let $D$ be a $G$-graded division algebra with support $T$ and $V$ a $G$-graded right $D$-module over an algebraically closed field $\FF$ satisfying the conditions of Proposition \ref{p:GDivAlgWithInv}. Let $\varphi_0$ be an involution of $D$ such that $(D,\varphi_0)$ is a simple algebra with involution and $B\colon V\times V\to D$ a nondegenerate hermitian or skew-hermitian $\varphi_{0}$-sesquilinear form such that $B\colon ^{\tau}V\times ^{\tau}V\to D$ is homogeneous of degree $(\subo,g)$. As shown in  Proposition \ref{p:GDivAlgWithInv} and Proposition \ref{p:RemoveTau} there is a homogeneous element $d\in D$ for which, after changing $\varphi_0$ by $\Ind(d)\circ \varphi_0$ and $B$ by $dB$, and using Proposition \ref{p:involutions} we can assume that $D$ is either of the form $\Dinv(T,\beta)$ or of the form $\Dinv(T,\beta)^{\ex}_t$. We write $\Sgn(B)=1$ if it is hermitian and $\Sgn(B)=-1$ if it is skew-hermitian.  We can decompose $V$ into isotypic components as:
	
	\begin{multline}\label{eq:isotypicComps}
		V=V_1\oplus...\oplus V_{m_0}\oplus( V'_{m_0+1}\oplus V''_{m_0+1})\oplus...\oplus(V'_{k_0}\oplus V''_{k_0})\\
		\oplus W_1\oplus...\oplus W_{m_1}\oplus (W'_{m_1+1}\oplus W''_{m_1+1})\oplus...\oplus (W'_{k_1}\oplus W''_{k_1})
	\end{multline}
	in such way that each of the isotypic component $V_i$, $V'_i$ and $V''_i$ is isomorphic to a sum of submodules isomorphic to $D^{[(0,g)]}$ for a $g\in G$ which is fixed in each component and  each of the isotypic component $W_j$, $W'_j$ and $W''_j$ is isomorphic to a sum of submodules isomorphic to $D^{[(1,g)]}$ for a $g\in G$ which is fixed in each component, in such way that each $V_i$ and $W_j$ are self dual, and each of the pairs $V'_i$ and $V''_i$ or $W'_j$ and $W''_j$ are in  duality by $B$. This condition implies that $\dim_D(V'_i)=\dim_D(V''_i)$ and $\dim_D(W'_j)=\dim_D(W''_j)$ for each $i\in \{m_0+1,...,k_0\}$ and $j\in \{m_1+1,...,k_1\}$. We can order each $V_1,...,V_{m_0}$ such that the first $l_0$ have odd dimension and the remaining even dimension and similarly with each $W_1,...,W_{m_1}$.
	
	We write:
	
	\begin{equation}\label{eq:kappa0}
		\kappa_0=(q^{(0)}_1,...,q^{(0)}_{l_0},2q^{(0)}_{l_0+1},...,2q^{(0)}_{m_0},q_{m_0+1}^{(0)},q_{m_0+1}^{(0)},...,q^{(0)}_{k_0},q^{(0)}_{k_0})\in \ZZ_{> 0}^{k_0}
	\end{equation}
	to be the vector with the dimensions of each $V_i$, $V'_i$ and $V''_i$ and 
	
	\begin{equation}\label{eq:kappa1}
		\kappa_1=(q^{(1)}_1,...,q^{(1)}_{l_1},2q^{(1)}_{l_1+1},...,2q^{(1)}_{m_1},q_{m_1+1}^{(1)},q_{m_1+1}^{(1)},...,q^{(1)}_{k_1},q^{(1)}_{k_1})\in \ZZ_{> 0}^{k_1}
	\end{equation}
	to be the vector with the dimensions of each $W_i$, $W'_i$ and $W''_i$. Here $q^{(0)}_1,...,q^{(0)}_{l_0},$  and $q^{(1)}_1,...,q^{(1)}_{l_1}$ are odd numbers, and write:
	
	\begin{equation}\label{eq:gamma0}
		\gamma_0=(g^{(0)}_1,...,g^{(0)}_{l_0},g^{(0)}_{l_0+1},...,g^{(0)}_{m_0},g'^{(0)}_{m_0+1},g''^{(0)}_{m_0+1},...,g'^{(0)}_{k_0},g''^{(0)}_{k_0})\in G^{2k_0-m_0}
	\end{equation}
	and 
	
	\begin{equation}\label{eq:gamma1}
		\gamma_1=(g^{(1)}_1,...,g^{(1)}_{l_1},g^{(1)}_{l_1+1},...,g^{(1)}_{m_0},g'^{(1)}_{m_1+1},g''^{(1)}_{m_1+1},...,g'^{(1)}_{k_1},g''^{(1)}_{k_1})\in G^{2k_1-m_1}
	\end{equation}
	
	For group elements satisfying that $\supp V_i=g^{(0)}_i T$, $\supp V'_i=g'^{(0)}_i T$, $\supp V''_i=g''^{(0)}_iT$, $\supp W_j=g^{(1)}_j T$, $\supp W'_j=g'^{(1)}_j T$ and $\supp W''_j=g''^{(1)}_j$ for each possible choice of $i\in\{1,...,k_0\}$ and $j\in\{1,...,k_1\}$. Given $h_1,h_2\in G$, $i,j\in\{0,1\}$, $v\in V_{(i,h_1)}$ and $w\in V_{(j,h_2)}$, we have that $B(v,w)=0$ unless $i=j$ and $h_1h_2g\in T$. Therefore, there are $t^{(0)}_1,...,t^{(0)}_{k_0},t^{(1)}_1,...,t^{(1)}_{k_1}\in T$ satisfying:
	
	\begin{multline*}
		(g_1^{(0)})^2t^{(0)}_1=...=(g_{m_0}^{(0)})^2t^{(0)}_{m_0}=g'^{(0)}_{m_0+1}g''^{(0)}_{m_0+1}t^{(0)}_{m_0+1}=...=g'^{(0)}_{k_0}g''^{(0)}_{k_0}t^{(0)}_{k_0}\\
		=	(g_1^{(1)})^2t^{(1)}_1=...=(g_{m_1}^{(1)})^2t^{(1)}_{m_1}=g'^{(1)}_{m_1+1}g''^{(1)}_{m_1+1}t^{(1)}_{m_1+1}=...=g'^{(1)}_{k_1}g''^{(1)}_{k_1}t^{(1)}_{k_1}=g^{-1}.
	\end{multline*} 
	
	If we change $g''^{(0)}_i$ with $g''^{(0)}_it^{(0)}_i$ for every $i\in\{m_0+1,...,k_0\}$ and  $g''^{(1)}_j$ with  $g''^{(1)}_it^{(1)}_i$ for every $i\in\{m_1+1,...,k_1\}$, we can assume that $\gamma_0$ and $\gamma_1$ satisfy:
	
	\begin{multline}\label{eq:degreesre1}
		(g_1^{(0)})^2t^{(0)}_1=...=(g_{m_0}^{(0)})^2t^{(0)}_{m_0}=g'^{(0)}_{m_0+1}g''^{(0)}_{m_0+1}=...=g'^{(0)}_{k_0}g''^{(0)}_{k_0}\\
		=	(g_1^{(1)})^2t^{(1)}_1=...=(g_{m_1}^{(1)})^2t^{(1)}_{m_1}=g'^{(1)}_{m_1+1}g''^{(1)}_{m_1+1}=...=g'^{(1)}_{k_1}g''^{(1)}_{k_1}=g^{-1}.
	\end{multline} 
	
	\begin{remark}
		With the previous notation, the elements \[t^{(0)}_1,...,t^{(0)}_{m_0},t^{(1)}_1,...,t^{(1)}_{m_1}\] and \[g''^{(0)}_{m_0+1},...,g''^{(0)}_{k_0},g''^{(1)}_{m_1+1},...,g''^{(1)}_{k_1},\] are determined by \[g^{(0)}_1,...,g^{(0)}_{m_0},g'^{(0)}_{m_0+1},...,g'^{(0)}_{k_0},g^{(1)}_1,...,g^{(1)}_{m_1},g'^{(1)}_{m_1+1},...,g'^{(1)}_{k_1},\] and $g$.
		
	\end{remark}

	\begin{df}
		Let $T'$ be an elementary $2$-subgroup of $\ZZ\times G$ and $\beta\colon T'\times T'\to \FF$ a nondegenerate alternating bicharacter. Let $\kappa_0$, $\kappa_1$, $\gamma_0$, $\gamma_1$ and $g\in G$ be as in \eqref{eq:kappa0}, \eqref{eq:kappa1}, \eqref{eq:gamma0} and \eqref{eq:gamma1} satisfying \eqref{eq:degreesre1} and let $\delta=\pm 1$. We denote $T=T'$ and define $\cM(G,T,\beta,\kappa_0,\kappa_1,\gamma_0, \gamma_1,\delta,g)$ as the graded algebra $\cM(G,\Dinv(T,\beta),\kappa_0,\kappa_1,\gamma_0, \gamma_1)$ with the involution given by $\varphi(X)=\Phi^{-1}X^t\Phi$ where $\Phi$ is the matrix given in the following block-diagonal form:
		
		\begin{multline}\label{eq:MatrixInv1}
			\Phi=\sum_{i=1}^{l_0} \mathrm{I}_{q^{(0)}_{i}} \otimes X_{t^{(0)}_i}\oplus \sum_{i=l_0+1}^{m_0} \mathrm{S^{(0)}}_{i} \otimes X_{t^{(0)}_i}\oplus\sum _{m_0+1}^{k_0}\begin{pmatrix}
				0& I_{q^{(0)}_i}\\\delta I_{q^{(0)}_i}&0\\
			\end{pmatrix}\\ \oplus \sum_{i=1}^{l_1} \mathrm{I}_{q^{(1)}_{i}} \otimes X_{t^{(1)}_i}\oplus \sum_{i=l_1+1}^{m_1} \mathrm{S^{(1)}}_{i} \otimes X_{t^{(1)}_i}\oplus\sum _{m_1+1}^{k_1}\begin{pmatrix}
				0& I_{q^{(1)}_i}\\\delta I_{q^{(1)}_i}&0\\
			\end{pmatrix}
		\end{multline}
		where for $i\in \{l_0+1,...,m_0\}$, $S^{(0)}_i=\mathrm{I}_{2q_{i}^{(0)}}$ or $S^{(0)}_i=\begin{pmatrix} 0&\mathrm{I}_{q_{i}^{(0)}}\\-\mathrm{I}_{q_{i}^{(0)}}&0
		\end{pmatrix}$, we write $\sgn(S^{(0)}_i)=1$ in the first case and $\sgn(S^{(0)}_i)=-1$ in the second case, and similarly, for  $i\in \{l_1+1,...,m_1\}$, $S^{(1)}_i=\mathrm{I}_{2q_{i}^{(1)}}$ or $S^{(1)}_i=\begin{pmatrix} 0&\mathrm{I}_{q_{i}^{(1)}}\\-\mathrm{I}_{q_{i}^{(1)}}&0
		\end{pmatrix}$, and we write $\sgn(S^{(1)}_i)=1$ in the first case and $\sgn(S^{(1)}_i)=-1$ in the second case. All this subject to:
		
		\begin{multline}\label{eq:RestrictionsInv1}
			\delta=\beta(t^{(0)}_1)=...=\beta(t^{(0)}_{l_0})=\sgn(S^{(0)}_{l_0+1})\beta(t^{(0)}_{l_0+1})=...=\sgn(S^{(0)}_{m_0})\beta(t^{(0)}_{m_0})\\ =\beta(t^{(1)}_1)=...=\beta(t^{(1)}_{l_1})=\sgn(S^{(1)}_{l_1+1})\beta(t^{(1)}_{l_1+1})=...=\sgn(S^{(1)}_{m_1})\beta(t^{(1)}_{m_1})
		\end{multline}
		(recall that $\beta$ is the quadratic form defined as in Definition \ref{df:Dinv}). 
		
		With the same notation, for an order $2$ element $t\in G\setminus T'$, we denote $T=T'\langle t\rangle$ and define $\cM(G,T',t,\beta,\kappa_0,\kappa_1,\gamma_0, \gamma_1,\delta,g)$ to be $\cM(G,\Dinv(T',\beta)^{\ex}_t,\kappa_0,\kappa_1,\gamma_0, \gamma_1)$ with the involution given by the expression $\varphi(X)=\Phi^{-1}\ex(X)^t\Phi$ where $\ex$ is applied componentwise and $\Phi$ is the matrix given in the following block-diagonal form:
		
		\begin{multline}\label{eq:MatrixInv2}
			\Phi=\sum_{i=1}^{l_0} \mathrm{I}_{q^{(0)}_{i}} \otimes Y_{t^{(0)}_i}\oplus \sum_{i=l_0+1}^{m_0} \mathrm{S^{(0)}}_{i} \otimes Y_{t^{(0)}_i}\oplus\sum _{m_0+1}^{k_0}\begin{pmatrix}
				0& I_{q^{(0)}_i}\\\delta I_{q^{(0)}_i}&0\\
			\end{pmatrix}\\ \oplus \sum_{i=1}^{l_1} \mathrm{I}_{q^{(1)}_{i}} \otimes Y_{t^{(1)}_i}\oplus \sum_{i=l_1+1}^{m_1} \mathrm{S^{(1)}}_{i} \otimes Y_{t^{(1)}_i}\oplus\sum _{m_1+1}^{k_1}\begin{pmatrix}
				0& I_{q^{(1)}_i}\\\delta I_{q^{(1)}_i}&0\\
			\end{pmatrix}
		\end{multline}
		
		Where for $i\in \{l_0+1,...,m_0\}$, $S^{(0)}_i=\mathrm{I}_{2q_{i}^{(0)}}$ or $S^{(0)}_i=\begin{pmatrix} 0&\mathrm{I}_{q_{i}^{(0)}}\\-\mathrm{I}_{q_{i}^{(0)}}&0
		\end{pmatrix}$, we write $\sgn(S^{(0)}_i)=1$ in the first case and $\sgn(S^{(0)}_i)=-1$ in the second case, and similarly, for  $i\in \{l_1+1,...,m_1\}$, $S^{(1)}_i=\mathrm{I}_{2q_{i}^{(1)}}$ or $S^{(1)}_i=\begin{pmatrix} 0&\mathrm{I}_{q_{i}^{(1)}}\\-\mathrm{I}_{q_{i}^{(1)}}&0
		\end{pmatrix}$, and we write $\sgn(S^{(1)}_i)=1$ in the first case and $\sgn(S^{(1)}_i)=-1$ in the second case. All this subject to:
		
		\begin{multline}\label{eq:RestrictionsInv2}
			\delta=\beta^{[t]}(t^{(0)}_1)=...=\beta^{[t]}(t^{(0)}_{l_0})=\sgn(S^{(0)}_{l_0+1})\beta^{[t]}(t^{(0)}_{l_0+1})=...=\sgn(S^{(0)}_{m_0})\beta^{[t]}(t^{(0)}_{m_0})\\ =\beta^{[t]}(t^{(1)}_1)=...=\beta^{[t]}(t^{(1)}_{l_1})=\sgn(S^{(1)}_{l_1+1})\beta^{[t]}(t^{(1)}_{l_1+1})=...=\sgn(S^{(1)}_{m_1})\beta^{[t]}(t^{(1)}_{m_1}).
		\end{multline} 
	\end{df}
	
	\begin{proposition}\label{p:GradedSimple}
		Let $D$ be a $G$-graded division algebra with support $T$ and $V$ a $G$-graded right $D$-module over an algebraically closed field $\FF$ satisfying the conditions of Proposition \ref{p:GDivAlgWithInv}.  Let $\varphi_0$ be an involution of $D$ such that $(D,\varphi_0)$ is a simple algebra with involution and $B\colon V\times V\to D$ a nondegenerate hermitian or skew-hermitian $\varphi_{0}$-sesquilinear form such that $B\colon ^{\tau}V\times ^{\tau}V\to D$ is homogeneous of degree $(\subo,g)$ for $\tau$ as in $\eqref{eq:tau}$.
		
		\begin{itemize}
			\item[(1)]If $(D,\varphi_0)$ is isomorphic to $\Dinv(T,\beta)$, then $\End(G,D,V,\varphi_0,B)$ is isomorphic to $\cM(G,T,\beta,\kappa_0,\kappa_1,\gamma_0, \gamma_1,\delta,g)$ with $\kappa_0$, $\kappa_1$, $\gamma_0$, $\gamma_1$ and $g\in G$  as in \eqref{eq:kappa0}, \eqref{eq:kappa1}, \eqref{eq:gamma0} and \eqref{eq:gamma1} and $\delta=\sgn(B)$.
			
			\item[(2)]If $(D,\varphi_0)$ is isomorphic to $\Dinv(T',\beta)^{\ex}_t$ for some subgroup $T'$, then the graded algebra with involution $\End(G,D,V,\varphi_0,B)$ is isomorphic to the graded algebra with involution $\cM(G,T,t,\beta,\kappa_0,\kappa_1,\gamma_0, \gamma_1,\delta,g)$ with $\kappa_0$, $\kappa_1$, $\gamma_0$, $\gamma_1$ and $g\in G$  as in \eqref{eq:kappa0}, \eqref{eq:kappa1}, \eqref{eq:gamma0} and \eqref{eq:gamma1} and $\delta=\sgn(B)$. Moreover, with the same choice of $D,V$ and $\varphi_0$, we can choose $B$ so that $\sgn(B)=1$. 
			
		\end{itemize}
	\end{proposition}
	
	\begin{proof}
		We will abuse of the notation and in case $(D,\varphi_0)$ is isomorphic to $\Dinv(T',\beta)^{\ex}_t$, we will  denote $\beta^{[t]}$ just by $\beta$  in order to work with $(1)$ and $(2)$ at the same time. Also, for every $s\in T$, we will denote $Z_s=X_s$ whenever $(D,\varphi_0)$ is isomorphic to $\Dinv(T,\beta)$, and $Z_s=Y_s$ whenever it is isomorphic to $\Dinv(T',\beta)^{\ex}_t$. As shown in \cite[(2.19)]{EKmon}, if we take a homogeneous basis $\{v_1,...,v_r\}$ of $V$, the involution $\varphi$ of $\End(G,D,V,\varphi_0,B)$, can be written in matrix form as $\varphi(X)=\Phi^{-1}\varphi_0(X)^t\Phi$, where the coordinate $i,j$ of $\Phi$ is $B(v_i,v_j)$. Hence, we are going to find a good homogeneous basis to prove the proposition.
		
		For each $i\in\{1,...,m_0\}$ and $j\in \{1,...,m_1\}$, denote $\overline{V}_i=V_{g^{(0)}_i}$ and $\overline{W}_j=V_{g^{(1)}_j}$. In that way, we can identify $V_i$ with $\overline{V_i}\otimes D$ and $W_j$ with $\overline{W}_j\otimes D$. Then, we have:
		
		\begin{equation*}
			\begin{split}
				B(u,v)=B^{(0)}_i(u,v)Z_{t^{(0)}_i}& \text{ for all }u,v\in \overline{V_i}\\
				B(u,v)=B^{(1)}_j(u,v)Z_{t^{(1)}_i}& \text{ for all }u,v\in \overline{W_j}.
			\end{split}
		\end{equation*}

		For some nondegenerate bilinear forms $B^{(0)}_i\colon \overline{V_i}\times \overline{V_i}\to \FF$ and $B^{(1)}_i\colon \overline{W_i}\times \overline{W_i}\to \FF$. In case $\sgn(B)=\beta(t^{(0)}_i)$ the fact that $B$ is hermitian or skew-hermitian implies that $B^{(0)}_i$ is symmetric. Thus, we can take a basis such that the matrix of  $B^{(0)}_i$ with respect to this basis is $\mathrm{I}_{\dim_{\FF}(\overline{V}_i)}$, otherwise, if $\sgn(B)=-\beta(t^{(0)}_i)$, it implies that $B^{(0)}_i$ is skew-symmetric, in which case, the dimension of $\overline{V}_i$ is even, and we denote $2k_i=\dim_{\FF}(\overline{V}_i)$. In these conditions, we can take a basis $\overline{V_i}$ satisfying that the matrix of $B^{(0)}_i$ is $\begin{pmatrix} 0&\mathrm{I}_{k_i}\\-\mathrm{I}_{k_i}&0
		\end{pmatrix}$. We proceed in the same way for each $B^{(1)}_j$.
		
		For every $i\in \{m_0+1,...,k_0\}$ and $j\in \{m_1+1,...,k_1\}$, we denote $\overline{V}'_i=V_{g'^{(0)}_i}$,  $\overline{V}''_i=V_{g''^{(0)}_i}$, $\overline{W}'_j=V_{g'^{(1)}_j}$ and $\overline{W}''_j=V_{g''^{(1)}_j}$. In that way, we can identify $V'_i$ with $\overline{V}'_i\otimes D$, $V''_i$ with $\overline{V}''_i\otimes D$, $w'_j$ with $\overline{W}'_j\otimes D$ and $V''_i$ with $\overline{W}''_i\otimes D$. And we can show that there are nondegenerate pairings $B^{(0)}_i\colon \overline{V}'_i\times \overline{V}''_i\to \FF$ and $B^{(1)}_i\colon \overline{W}'_j\times \overline{W}''_j\to \FF$  such that:
		
		\begin{equation*}
			\begin{split}
				B(u,v)=B^{(0)}_i(u,v)Z_{t^{(0)}_i}& \text{ for all }u\in \overline{V}_1 \text{ and } v\in \overline{V}''_i\\
				B(u,v)=B^{(1)}_j(u,v)Z_{t^{(1)}_i}& \text{ for all }u\in \overline{W}'_j\text{ and }v\in \overline{W}''_j.
			\end{split}
		\end{equation*}
		
		Thus, we can take a basis of $\overline{V}'_i$ and a dual basis of $\overline{V}''_i$ with respect to $B_i^{(0)}$ and also a basis for $\overline{W}'_i$ and a dual basis of $\overline{W}''_i$ with respect to $B_i^{(0)}$. Ordering the basis of each subspace as in \eqref{eq:isotypicComps}, we can write the algebra in matrix form and get the result.

		Finally, if $(D,\varphi_0)$ is isomorphic to $\Dinv(T',\beta)^{\ex}_t$ and $\sgn(B)=-1$, due to Proposition \ref{p:involutions} \[\End(G,\cD,V,\varphi_0,B)=\End(G,\cD,V,\mathrm{Int}((1,-1))\circ \varphi_0,(1,-1)B)\] but $\mathrm{Int}(1,-1)\circ \varphi_0=\varphi_0$ and $\sgn((1,-1)B)=1$.
		
	\end{proof}

	\begin{theorem}\label{t:ClassificationAssociativeCase}
		Let $(A,\varphi,\Gamma)$ be an $\Omega$-algebra satisfying $(T1)-(T4)$ over an algebraically closed field $\FF$. Then, one of the following holds:
		
		\begin{itemize}
			\item[(1)] There is an elementary $2$-subgroup $T$ of $G$, a nondegenerate alternating bicharacter $\beta\colon T\times T\to \FF$, $k_0,k_1>0$, $\gamma_0\in G^{k_0}$ and $\gamma_1\in G^{k_1}$ consisting on distinct elements modulo $T$,  $\kappa_0\in \ZZ^{k_0}_{>0}$ and $\kappa_1\in \ZZ^{k_1}_{>0}$ such that $(A,\varphi,\Gamma)\cong \cM(G,D(T,\beta),\kappa_0,\kappa_1,\gamma_0,\gamma_1)^{\ex}$.
			\item[(2)] There are, an elementary $2$-subgroup $T$ of $G$, a nondegenerate alternating bicharacter $\beta\colon T\times T\to \FF$,  $\kappa_0$, $\kappa_1$, $\gamma_0$, $\gamma_1$ and $g\in G$  as in \eqref{eq:kappa0}, \eqref{eq:kappa1}, \eqref{eq:gamma0} and \eqref{eq:gamma1} with $\kappa_0,\kappa_1\neq (0)$ and $\delta=\pm 1$ such that we have the isomorphism $(A,\varphi,\Gamma)\cong \cM(G,T,\beta,\kappa_0,\kappa_1,\gamma_0, \gamma_1,\delta,g)$.
			
			\item[(3)] There are, an elementary $2$-subgroup $T$ of $G$, a nondegenerate alternating bicharacter $\beta\colon T\times T\to \FF$, an order $2$ element $t\in G\setminus T$,  $\kappa_0$, $\kappa_1$, $\gamma_0$, $\gamma_1$ and $g\in G$  as in \eqref{eq:kappa0}, \eqref{eq:kappa1}, \eqref{eq:gamma0} and \eqref{eq:gamma1}  with $\kappa_0,\kappa_1\neq (0)$  such that $(A,\varphi,\Gamma)\cong \cM(G,T,\beta,t,\kappa_0,\kappa_1,\gamma_0, \gamma_1,1,g)$.
		\end{itemize}
		
		Moreover, the algebras on each item are non-isomorphic and:
		
		\begin{itemize}
			\item[(i)] $\cM(G,D(T,\beta),\kappa_0,\kappa_1,\gamma_0,\gamma_1)^{\ex}\cong \cM(G,D(T',\beta'),\kappa'_0,\kappa'_1,\gamma'_0,\gamma'_1)^{\ex}$ if and only if  either:
			\begin{itemize}
				\item[(a)] $T=T'$, $\beta=\beta'$ and there is $g\in G$ such that $\Xi(\kappa_0,\gamma_0)=g\Xi(\kappa'_0,\gamma'_0)$ and $\Xi(\kappa_1,\gamma_1)=g\Xi(\kappa'_1,\gamma'_1)$ 
				\item[(b)] $T=T'$, $\beta=\beta'$ and there is $g\in G$ such that $\Xi(\kappa_0,\gamma_0)=g\Xi(\kappa'_0,\gamma'^{-1}_0)$ and $\Xi(\kappa_1,\gamma_1)=g\Xi(\kappa'^{-1}_1,\gamma'_1)$. 
			\end{itemize}
			\item[(ii)] $\cM(G,T,\beta,\kappa_0,\kappa_1,\gamma_0, \gamma_1,\delta,g)\cong \cM(G,T',\beta',\kappa_0',\kappa_1',\gamma_0', \gamma_1',\delta',g')$ if and only if $T=T'$, $\beta=\beta'$, $\delta=\delta'$ and there is $g''\in G$ such that $\Xi(\kappa_0,\gamma_0)=g''\Xi(\kappa'_0,\gamma'_0)$, $\Xi(\kappa_1,\gamma_1)=g''\Xi(\kappa'_1,\gamma'_1)$ and $g=g'g''^{-1}$.
			
			\item[(iii)] $\cM(G,T,\beta,t,\kappa_0,\kappa_1,\gamma_0, \gamma_1,1,g)\cong \cM(G,T',\beta',t',\kappa'_0,\kappa'_1,\gamma'_0, \gamma'_1,1,g')$ if and only if $T\langle t\rangle=T'\langle t'\rangle$, $t=t'$, $\beta^{[t]}=\beta^{[t']}$ and there is $g''\in G$ such that $\Xi(\kappa_0,\gamma_0)=g''\Xi(\kappa'_0,\gamma'_0)$, $\Xi(\kappa_1,\gamma_1)=g''\Xi(\kappa'_1,\gamma'_1)$ and $g=g'g''^{-1}$.
		\end{itemize}
	\end{theorem}

	\begin{proof}
		In order to  prove the first part, assume that $(A,\varphi,\Gamma)$ is an $\Omega$-algebra satisfying $(T1)-(T4)$ over an algebraically closed field $\FF$. In case $(A,\varphi,\Gamma)$ is not graded simple, Theorems \ref{t:NonGradedSimple} and \ref{t:DivAlgClass} imply $(1)$. In case $A$ is simple, due to Proposition \ref{p:GradedSimple}, we get $(2)$ and $(3)$.
		
		We are left with the second part. In order to prove it, we first notice that the algebras on each item are not isomorphic since in $(1)$, $(A,\Gamma)$ is not graded-simple, in $(2)$, $A$ is simple, hence $(A,\Gamma)$ is graded-simple and in $(3)$ $(A,\Gamma)$ is graded-simple but $A$ is not simple. Hence, let's prove item by item:
		
		\begin{itemize}
			\item[(i)]This is a consequence again of Theorems \ref{t:NonGradedSimple} and \ref{t:DivAlgClass} together with Lemma \ref{l:DivAlgOp} and Remark \ref{rm:bicharacterOp}.
			
			\item[(ii)] Let $D$, $V$, $\varphi_0$ and $B$ be such that the isomorphism $\End(G,D,V,\varphi_0,B)\cong\cM(G,T,\beta,\kappa_0,\kappa_1,\gamma_0, \gamma_1,\delta,g)$ holds. If $T=T'$, $\beta=\beta'$, $\delta=\delta'$ and there is $g''\in G$ such that $\Xi(\kappa_0,\gamma_0)=g\Xi(\kappa'_0,\gamma'_0)$, $\Xi(\kappa_1,\gamma_1)=g\Xi(\kappa'_1,\gamma'_1)$ and $g=g'g''^{-1}$. Then, by the construction, shown before, it follows \[\End(G,D,V^{[(0,g'')]},\varphi_0,B)\cong\cM(G,T',\beta',\kappa'_0,\kappa'_1,\gamma'_0, \gamma'_1,\delta',g')\] and it also follows, that we have the isomorphism \[\End(G,D,V,\varphi_0,B)\cong \End(G,D,V^{[(0,g'')]},\varphi_0,B)\] via the isomorphism induced by $(\id_{D},\id_V)$. The converse is just a consequence of \cite[Theorem 2.64]{EKmon} since here we are considering a subgroup of the automorphisms in the theorem mentioned before.
			\item[(iii)] If $T\langle t\rangle=T'\langle t'\rangle$, $t=t'$, $\beta^{[t]}=\beta^{[t']}$ and there is $g''\in G$ such that $\Xi(\kappa_0,\gamma_0)=g\Xi(\kappa'_0,\gamma'_0)$, $\Xi(\kappa_1,\gamma_1)=g\Xi(\kappa'_1,\gamma'_1)$ and $g=g'g''^{-1}$, the prove is as before.

			For the converse, Let $D$, $V$, $\varphi_0$ and $B$ be such that $\End(G,D,V,\varphi_0,B)\cong\cM(G,T,\beta,\kappa_0,\kappa_1,\gamma_0, \gamma_1,1,g)$. Due to \cite[Theorem 2.10]{EKmon}, the pair $(D,V)$ is determined up to isomorphism by the grading that we have on the $\Omega$-algebra
			$\cM(G,T,\beta,\kappa_0,\kappa_1,\gamma_0, \gamma_1,1,g)$ and a shift on the grading of $V$ by an element $(0,g'')\in \ZZ\times G$.  Due to Proposition \ref{p:GradedSimple} the sesquilinear form $B$ is determided up to one homogeneous element of the center of the algebra $\End_{D}(V)$ by the involution of $\End(G,D,V,\varphi_0,B)$ and the choice of $\varphi_0$. The restriction $\sgn(B)=1$, implies that it is uniquely determined by scalar. Therefore, the decomposition \eqref{eq:isotypicComps} is determined up to a permutation of $V_1,...,V_{m_0}$, a permutation of $V'_{m_0+1}, V''_{m_0+1},..., V'_{k_0},V''_{k_0}$ preserving the pairing of the $V'_i$ and $V''_i$, a permutation of $W_1,...,W_{m_1}$, and a permutation of $W'_{m_1+1}, W''_{m_1+1},...$  $, W'_{k_1},W''_{k_1}$ preserving the pairing of the $V'_i$ and $V''_i$. If we shift the grading on $V$ by $(0,g'')$, the degree of $^{[\tau]}B$ changes to $(\subo,gg''^{-1})$. Thus, the elements the elements $t^{(0)}_1,...,t^{(0)}_{m_0},t^{(1)}_1,...,t^{(1)}_{m_1}$ are uniquely determined by $B$ since the equation \eqref{eq:degreesre1} is invariant under these shifts and any changes of the elements $g^{(0)}_1,...,g^{(0)}_{m_0},g^{(1)}_1,...,g^{(1)}_{m_1}$ by other elements on their cosets (since $T$ is an elementary $2$-group). Therefore, this argument together with Theorem \ref{t:NonSimpleGradedDivAlg} implies that  in case \[\cM(G,T,\beta,t,\kappa_0,\kappa_1,\gamma_0, \gamma_1,1,g')\cong \cM(G,T',\beta',t',\kappa'_0,\kappa'_1,\gamma'_0, \gamma'_1,1,g')\], then $T\langle t\rangle=T'\langle t'\rangle$, $t=t'$, $\beta^{[t]}=\beta^{[t']}$ and there is $g''\in G$ such that $\Xi(\kappa_0,\gamma_0)=g\Xi(\kappa'_0,\gamma'_0)$, $\Xi(\kappa_1,\gamma_1)=g\Xi(\kappa'_1,\gamma'_1)$ and $g=g'g''^{-1}$.
			
		\end{itemize}

	\end{proof}
	
	\section{Gradings on simple AT2 over algebraically closed fields}\label{sec:GradingsAt2}
	
	If $(A,\varphi,\Gamma)$ is a triple satisfying $(T1)-(T4)$, and we denote $^{\pi_1}\Gamma\colon A=A_{-1}\oplus A_0\oplus A_1$ by $\bold{GrW}(A,\varphi,\Gamma)$ the graded associative triple $(W,\{\cdots\},\theta^{-1}(^{\pi_2}\Gamma))$, where $(W,\{\cdots\})$ is the triple system $\cW(A)$, $\theta$ is the isomorphism from Theorem \ref{t:AutSchemeIso} and $\theta^{-1}(^{\pi_2}\Gamma)$ is  the grading from Proposition \ref{p:TransferThm}.
	
	Given finitely generated abelian group $G$ we will define the following triple systems:
	
	\begin{itemize}
		\item[(1)] For an elementary $2$-subgroup $T$ of $G$, a nondegenerate alternating bicharacter $\beta\colon T\times T\to \FF$, $k_0,k_1>0$, $\gamma_0\in G^{k_0}$ and $\gamma_1\in G^{k_1}$ consisting on distinct elements modulo $T$,  $\kappa_0\in \ZZ^{k_0}_{>0}$ and $\kappa_1\in \ZZ^{k_1}_{>0}$ we denote the graded triple system $\textbf{GrW}(\cM(G,D(T,\beta),\kappa_0,\kappa_1,\gamma_0,\gamma_1)^{\ex})$ by \[\textbf{GrW}(G,D(T,\beta),\kappa_0,\kappa_1,\gamma_0,\gamma_1)^{\ex}.\]
		\item[(2)] For an elementary $2$-subgroup $T$ of $G$, a nondegenerate alternating bicharacter $\beta\colon T\times T\to \FF$,  $\kappa_0$, $\kappa_1$, $\gamma_0$, $\gamma_1$ and $g\in G$  as in \eqref{eq:kappa0}, \eqref{eq:kappa1}, \eqref{eq:gamma0} and \eqref{eq:gamma1} with $\kappa_0,\kappa_1\neq (0)$ and $\delta=\pm 1$ we denote the graded triple system $\textbf{GrW}(\cM(G,T,\beta,\kappa_0,\kappa_1,\gamma_0, \gamma_1,\delta,g))$ by: \[\textbf{GrW}(G,T,\beta,\kappa_0,\kappa_1,\gamma_0, \gamma_1,\delta,g)\]
		
		\item[(3)] For an elementary $2$-subgroup $T$ of $G$, a nondegenerate alternating bicharacter $\beta\colon T\times T\to \FF$, an order $2$ element $t\in G\setminus T$,  $\kappa_0$, $\kappa_1$, $\gamma_0$, $\gamma_1$ and $g\in G$  as in \eqref{eq:kappa0}, \eqref{eq:kappa1}, \eqref{eq:gamma0} and \eqref{eq:gamma1}  with $\kappa_0,\kappa_1\neq (0)$  we denote the graded triple system  $\textbf{GrW}(\cM(G,T,\beta,t,\kappa_0,\kappa_1,\gamma_0, \gamma_1,1,g))$ by \[\textbf{GrW}(G,T,\beta,t,\kappa_0,\kappa_1,\gamma_0, \gamma_1,g)\]
	\end{itemize}

	\begin{example}\label{ex:GrATtwo}
		With the notation of Example \ref{ex:ATtwoex}  $\bold{GrW}(A,\varphi,\Gamma)$ is 
		$\cW(A,\varphi,\Gamma)$ with the grading given by: $\deg(X)=g$ if:
		
		\[\begin{pmatrix}
			0&X\\0&0\end{pmatrix}\in A_g\]
	\end{example}
	
	\begin{example}\label{ex:GrATtwoex}
		With the notation of Example \ref{ex:ATtwoex}  $\bold{GrW}(A,\varphi,\Gamma)$ is 
		$\cW(A,\varphi,\Gamma)$ with the grading given by: $\deg(X_1,X_2)=g$ if:
		
		\[\left(\begin{pmatrix}
			0&X_1\\0&0\end{pmatrix},\begin{pmatrix}
			0&X_2\\0&0\end{pmatrix}\right)\in A_g\]
	\end{example}

	\begin{theorem}\label{t:ClassificationTripleSystem}
		Let $(W,\{\cdots\},\Gamma)$ be an simple triple system with a $G$-grading $\Gamma$ over an algebraically closed field $\FF$. Then, one of the following holds:
		
		\begin{itemize}
			\item[(1)] There is an elementary $2$-subgroup $T$ of $G$, a nondegenerate alternating bicharacter $\beta\colon T\times T\to \FF$, $k_0,k_1>0$, $\gamma_0\in G^{k_0}$ and $\gamma_1\in G^{k_1}$ consisting on distinct elements modulo $T$,  $\kappa_0\in \ZZ^{k_0}_{>0}$ and $\kappa_1\in \ZZ^{k_1}_{>0}$ such that $(W,\{\cdots\},\Gamma)\cong \textbf{GrW}(G,D(T,\beta),\kappa_0,\kappa_1,\gamma_0,\gamma_1)^{\ex}$.
			\item[(2)] There are, an elementary $2$-subgroup $T$ of $G$, a nondegenerate alternating bicharacter $\beta\colon T\times T\to \FF$,  $\kappa_0$, $\kappa_1$, $\gamma_0$, $\gamma_1$ and $g\in G$  as in \eqref{eq:kappa0}, \eqref{eq:kappa1}, \eqref{eq:gamma0} and \eqref{eq:gamma1} with $\kappa_0,\kappa_1\neq (0)$ and $\delta=\pm 1$ such that we have the isomorphism $(W,\{\cdots\},\Gamma)\cong \textbf{GrW}(G,T,\beta,\kappa_0,\kappa_1,\gamma_0, \gamma_1,\delta,g)$.
			
			\item[(3)] There are, an elementary $2$-subgroup $T$ of $G$, a nondegenerate alternating bicharacter $\beta\colon T\times T\to \FF$, an order $2$ element $t\in G\setminus T$,  $\kappa_0$, $\kappa_1$, $\gamma_0$, $\gamma_1$ and $g\in G$  as in \eqref{eq:kappa0}, \eqref{eq:kappa1}, \eqref{eq:gamma0} and \eqref{eq:gamma1}  with $\kappa_0,\kappa_1\neq (0)$  such that $(W,\{\cdots\},\Gamma)\cong \textbf{GrW}(G,T,\beta,t,\kappa_0,\kappa_1,\gamma_0, \gamma_1,g)$.
		\end{itemize}
		
		Moreover, the algebras on each item are non isomorphic and:
		
		\begin{itemize}
			\item[(i)] $\textbf{GrW}(G,D(T,\beta),\kappa_0,\kappa_1,\gamma_0,\gamma_1)^{\ex}\cong \textbf{GrW}(G,D(T',\beta'),\kappa'_0,\kappa'_1,\gamma'_0,\gamma'_1)^{\ex}$ if and only if  either:
			\begin{itemize}
				\item[(a)] $T=T'$, $\beta=\beta'$ and there is $g\in G$ such that $\Xi(\kappa_0,\gamma_0)=g\Xi(\kappa'_0,\gamma'_0)$ and $\Xi(\kappa_1,\gamma_1)=g\Xi(\kappa'_1,\gamma'_1)$ 
				\item[(b)] $T=T'$, $\beta=\beta'$ and there is $g\in G$ such that $\Xi(\kappa_0,\gamma_0)=g\Xi(\kappa'_0,\gamma'^{-1}_0)$ and $\Xi(\kappa_1,\gamma_1)=g\Xi(\kappa'^{-1}_1,\gamma'_1)$. 
			\end{itemize}
			\item[(ii)] $\textbf{GrW}(G,T,\beta,\kappa_0,\kappa_1,\gamma_0, \gamma_1,\delta,g)\cong \textbf{GrW}(G,T',\beta',\kappa_0',\kappa_1',\gamma_0', \gamma_1',\delta',g')$ if and only if $T=T'$, $\beta=\beta'$, $\delta=\delta'$ and there is $g''\in G$ such that $\Xi(\kappa_0,\gamma_0)=g''\Xi(\kappa'_0,\gamma'_0)$, $\Xi(\kappa_1,\gamma_1)=g''\Xi(\kappa'_1,\gamma'_1)$ and $g=g'g''^{-1}$.
			
			\item[(iii)] $\textbf{GrW}(G,T,\beta,t,\kappa_0,\kappa_1,\gamma_0, \gamma_1,1,g)\cong \textbf{GrW}(G,T',\beta',t',\kappa'_0,\kappa'_1,\gamma'_0, \gamma'_1,1,g')$ if and only if $T\langle t\rangle=T'\langle t'\rangle$, $t=t'$, $\beta^{[t]}=\beta^{[t']}$ and there is $g''\in G$ such that $\Xi(\kappa_0,\gamma_0)=g''\Xi(\kappa'_0,\gamma'_0)$, $\Xi(\kappa_1,\gamma_1)=g''\Xi(\kappa'_1,\gamma'_1)$ and $g=g'g''^{-1}$.
		\end{itemize}
		\begin{proof}
			This result follows from Proposition \ref{p:TransferThm} and Theorems \ref{t:AutSchemeIso} and \ref{t:ClassificationAssociativeCase}
		\end{proof}
	\end{theorem}

	\section*{Acknowledgement}
	
	This work has been suported by the F.P.I. grant PRE2018-087018. It has been partially supported  by grant PID2021-123461NB-C21, funded by MCIN/AEI/-10.13039/501100011033 and by "ERDF A way of making Europe" and is also partially supported by Departamento de Ciencia, Universidad y Sociedad del Conocimiento del Gobierno de Arag{\'o}n (grant code: E22-23R: ``{\'A}lgebra y Geometr{\'i}a''). 
	
	This work is based on the results of my PhD thesis. I thank my supervisor Alberto Elduque for his useful suggestions.

\end{document}